\documentclass[12pt]{amsart}
\usepackage{amssymb}
\usepackage{amscd}

\newcommand{\Z}{\mathbb{Z}}
\newcommand{\N}{\mathbb{N}}
\newcommand{\C}{\mathbb{C}}
\newcommand{\Ss}{\mathbb{S}}
\newcommand{\Aa}{\mathbb{A}}

\newcommand{\Pp}{\mathbb{P}}

\newcommand{\B}{\mathbf{B}}

\newcommand{\Cmod}{\mathbb{C}\text{-}\mathbf{mod}}
\newcommand{\Var}{\mathbf{Var}}

\newcommand{\T}{\mathcal{T}}
\newcommand{\U}{\mathcal{U}}
\newcommand{\M}{\mathcal{M}}
\newcommand{\F}{\mathcal{F}}
\newcommand{\Mbar}{\overline{\mathcal{M}}}
\newcommand{\Mbarrn}{\overline{\mathcal{M}}(r,n)}
\newcommand{\CalO}{\mathcal{O}}
\newcommand{\CalA}{\mathcal{A}}

\newcommand{\CalP}{\mathcal{P}}

\newcommand{\Ind}{\mathrm{Ind}}
\newcommand{\End}{\mathrm{End}}
\newcommand{\Hom}{\mathrm{Hom}}
\newcommand{\noneq}{\natural}

\newcommand{\tr}{\mathrm{tr}}
\newcommand{\gtr}{\mathrm{tr}_{q,t}}
\newcommand{\gch}{\mathrm{ch}_{q,t}}
\newcommand{\ch}{\mathrm{ch}}
\newcommand{\Par}{\mathrm{Par}}
\newcommand{\Decomp}{\mathrm{Decomp}}
\newcommand{\Int}{\mathrm{Int}}
\newcommand{\Fibre}{\mathrm{Fibre}}
\newcommand{\Orb}{\mathrm{Orb}}
\newcommand{\isomto}{\overset{\sim}{\rightarrow}}
\newcommand{\strat}{\dashv}
\newtheorem{theorem}{Theorem}[section]
\newtheorem{corollary}[theorem]{Corollary}

\newtheorem{proposition}[theorem]{Proposition}
\numberwithin{equation}{section}

\theoremstyle{definition}

\begin{document}
\title[Representations of wreath products on cohomology]
{Representations of wreath products on cohomology of De~Concini-Procesi 
compactifications}
\author{Anthony Henderson}
\address{School of Mathematics and Statistics,
University of Sydney, NSW 2006, AUSTRALIA}
\email{anthonyh@maths.usyd.edu.au}
\thanks{This work was supported by Australian Research Council grant DP0344185}
%%%%%%%%%%%%%%%%%%%%%%%%%%%%%%%%%%%%%%%%%%%%%%%%%%
\begin{abstract}
The wreath product $W(r,n)$ of the cyclic group of order $r$
and the symmetric group $S_n$ acts on the corresponding projective
hyperplane complement, and on its wonderful compactification
as defined by De Concini and Procesi. We give a formula for the characters
of the representations of $W(r,n)$ on the cohomology groups of this
compactification,
extending the result of Ginzburg and Kapranov in the $r=1$ case.
As a corollary, we get a formula for the Betti numbers which generalizes
the result of Yuzvinsky in the $r=2$ case. Our method involves applying 
to the nested-set stratification
a generalization of Joyal's theory of tensor
species, which includes a link between polynomial functors and
plethysm for general $r$.
We also give a new proof of Lehrer's
formula for the representations of $W(r,n)$ on the cohomology groups
of the hyperplane complement.
\end{abstract}
%%%%%%%%%%%%%%%%%%%%%%%%%%%%%%%%%%%%%%%%%%%%%%%%%
\maketitle
%%%%%%%%%%%%%%%%%%%%%%%%%%%%%%%%%%%%%%%%%%%%%%%%
\section*{Introduction}
Let $r$ and $n$ be positive integers.
Write $W(r,n)$
for the \emph{wreath product} $\mu_r\wr S_n$ of the cyclic group $\mu_r$ of
complex $r$th roots of $1$ and the symmetric group $S_n$.
Let $V(r,n)$ be the standard vector space on which $W(r,n)$
acts irreducibly as a complex reflection group, with hyperplane
arrangement $\CalA(r,n)$ (in this context $W(r,n)$
is usually called
$G(r,1,n)$). Explicitly, $V(1,n)=\C^n/\C(1,1,\cdots,1)$,
$W(1,n)=S_n$ acts on $V(1,n)$ by permuting the coordinates, and
\[ \CalA(1,n)=\left\{\left.
\{(z_1,\cdots,z_n)\in\C^n\,|\,z_i=z_j\}/\C(1,1,\cdots,1)\,
\right| 1\leq i\neq j\leq n\right\}. \]
If $r\geq 2$, $V(r,n)=\C^n$, $W(r,n)$ acts by permutations of coordinates
composed with multiplying coordinates by $r$th roots of $1$, and
\[ \CalA(r,n)=\left\{\left.\{z_i=0\}\,\right| 1\leq i\leq n\right\}\cup
\left\{\left.\{z_i=\zeta z_j\}\,\right| 1\leq i\neq j\leq n,
\zeta\in\mu_r\right\}. \]
Write $\M(r,n)$ for the \emph{projective hyperplane complement}, i.e.\
the set of points in $\Pp(V(r,n))$ which (viewed as lines in $V(r,n)$)
lie in none of the hyperplanes in $\CalA(r,n)$; this is
a nonsingular irreducible affine complex variety of dimension $n-1$ (or $n-2$
if $r=1$, except that $\M(1,1)$ is empty). 
Clearly $W(r,n)$ acts on $\M(r,n)$. 
Taking singular (equivalently, de Rham)
cohomology, we obtain a graded representation
$H^\bullet(\M(r,n),\C)$ of $W(r,n)$, whose character has
been computed by Lehrer (see Theorems \ref{1lehrerthm} and
\ref{lehrerthm} below).

Note that $\M(1,n)$ for $n\geq 2$ can be identified with the moduli space
$\M_{0,n+1}$ of ordered configurations of $n+1$ points of $\Pp^1$.
This is because the last point can be taken to be the point at infinity,
so that what remains is an $n$-tuple of 
distinct complex numbers, modulo the simultaneous
action of the group of affine-linear transformations of $\C$. 
This point of view makes it clear that the action of $W(1,n)=S_n$
on $\M(1,n)$ can be extended to an action of $S_{n+1}$, but that
extension will not arise here.

In \cite[Section 4]{wonderful}, De Concini and Procesi have defined,
for any hyperplane arrangement $\CalA$ in a vector space $V$,
a ``wonderful'' compactification $\Mbar$ of the projective hyperplane 
complement $\M$;
this is a nonsingular irreducible projective variety containing $\M$
as an open subvariety, in which the complement of $\M$ is a divisor
with normal crossings. Actually the compactification comes in different 
versions, depending on the choice of a ``building set''. The version
we will use
is as follows. Let $L$ be the lattice of all intersections of hyperplanes
in $\CalA$, ordered by reverse inclusion.
Let $\F$ be the set of \emph{irreducible} elements of this lattice,
where $X\in L$ is irreducible if $X\neq V$ and the arrangement in
$V/X$ induced by the hyperplanes in $\CalA$ which contain $X$ is not
isomorphic to the product of two nontrivial arrangements. We assume
for simplicity that $\{0\}\in\F$, i.e.\ the arrangement $\CalA$ is essential
and irreducible. Collecting the obvious maps $\M\to\Pp(V/X)$ for all
$X\in\F$, we get a map
\[ \alpha:\M\to\prod_{X\in\F}\Pp(V/X), \]
which is an embedding since one of the factors on the right is
$\Pp(V)$. The compactification $\Mbar$ is defined to be
the closure of $\alpha(\M)$ in the projective variety
$\prod_{X\in\F}\Pp(V/X)$. In the case $\CalA=\CalA(1,n)$ for $n\geq 2$, 
the compactification
$\Mbar(1,n)$ is the usual compactification of $\M_{0,n+1}$ considered in 
algebraic geometry, namely the moduli space of stable 
genus $0$ curves with $n+1$ labelled points. For later convenience, we
redefine $\Mbar(1,1)$ to be a single point, not the empty set.

It is clear that whenever $\CalA$ is the hyperplane
arrangement of a complex reflection group $W$, the action of $W$ on $\M$
extends to an action on $\Mbar$. In particular, the
action of $W(r,n)$ on $\M(r,n)$ extends to its compactification $\Mbar(r,n)$.
The question then arises of computing the character of the graded
$W(r,n)$-representation $H^\bullet(\Mbar(r,n),\C)$.
This has been solved, to some extent, in the case $r=1$: there is a recursive
formula, essentially due to Ginzburg and Kapranov, allowing the character of
$H^\bullet(\overline{\M}(1,n),\C)$ to be computed from those for smaller $n$
(see Theorem \ref{ginzburgthm} below).
The main result of this paper (Theorem \ref{mainthm}) 
is a formula for the character of
$H^\bullet(\Mbarrn,\C)$, $r\geq 2$,
in terms of the characters of
$H^\bullet(\M(r,m),\C)$ and $H^\bullet(\overline{\M}(1,m),\C)$ for $m\leq n$.
As a corollary, we get a formula for the Betti numbers of
$\Mbar(r,n)$, which generalizes that found in the
$r=2$ case by Yuzvinsky (\cite[Theorem 5.1]{yuzvinsky}).

In \S1 we introduce some combinatorial notation which is needed to state
these formulas. Apart from some fairly natural generalizations of plethysm,
this is all standard wreath-product formalism in the spirit of
\cite[Chapter I, Appendix B]{macdonald}. In \S2 we state the formulas
for characters of representations of wreath products on cohomology which
are proved in this paper, and indicate their relationship with previous work.
To give the flavour of the main result, a sample calculation in the
$r=2$ case is given in \S3. 

Our derivation of this result relies on the nested-set stratification
of $\Mbar(r,n)$ described by De Concini and Procesi. 
Most of this paper is in fact a building-up of the correct
techniques for translating this stratification into a recursive formula for
cohomology; in order to keep track of the group actions, some
sophisticated combinatorial book-keeping is needed.
In \S5 we rephrase the
stratification so that the strata are indexed not by nested sets but by
a certain kind of trees with labelled vertices on which $\mu_r$ acts.
This brings the situation closer to the operadic viewpoint of
\cite{ginzburgkapranov}; but the standard theory of tensor species 
as developed by Joyal is no
longer sufficient to express the recursive nature of these trees in
the $r\geq 2$ case.

For this reason, in \S4 we generalize this 
theory to the case
of $r$-species and their linear analogues, which we call $\B_r$-modules.
As in the $r=1$ case there is a link with polynomial functors,
in this case from $\C\mu_r$-$\mathbf{mod}$ to $\C$-$\mathbf{mod}$.
This section revisits some of the concepts of Macdonald's paper
\cite{polyfunctors}, the crucial new ingredient
being the operation of precomposing
a polynomial functor $\C\mu_r$-$\mathbf{mod}\to\C$-$\mathbf{mod}$
with the functor $\C\mu_r$-$\mathbf{mod}\to\C\mu_r$-$\mathbf{mod}$
induced by a polynomial functor $\C$-$\mathbf{mod}\to\C$-$\mathbf{mod}$.
This operation, and the plethysm it corresponds to, are exactly what is
needed to translate the De Concini-Procesi stratification
into a statement about cohomology; but I believe that the results of \S4 
will be of independent interest.

In \S6, we present a jazzed-up version of this theory, where the
vector spaces have $(\N\times\N)$-gradings, and some tensor products
are not commutative but graded-commutative with respect to the
second $\N$-grading. In \S7, we apply this to the cohomology of our varieties,
where the second grading is by degree and the first is by half the Hodge 
weight. The fact that our strata are minimally pure in the sense of
\cite{dimcalehrer} allows us to take alternating sum with respect to degree,
and still distinguish the different cohomology groups by their Hodge weights.
Putting the pieces
of the argument together, we prove Theorems \ref{ginzburgthm} and
\ref{mainthm}. Finally, in \S8 we use this same technology to
give clean proofs of Lehrer's results for $\M(r,n)$ (this section
is based on the proof in the $r=1$ case given by Getzler).

\noindent
\textit{Acknowledgement. }
I am indebted to Gus Lehrer for suggesting this problem, and for
helpful conversations on these matters.
%%%%%%%%%%%%%%%%%%%%%%%%%%%%%%%%%%%%%%%%%%%%%%%%%%%%%%%
\tableofcontents
%%%%%%%%%%%%%%%%%%%%%%%%%%%%%%%%%%%%%%%%%%%%%%%%%%%%%%%
\section{Characteristics and Plethysm}
The combinatorial framework
for discussing representations of wreath products $G\wr S_n$ is given in
\cite[Chapter I, Appendix B]{macdonald}. Here we will only use the 
case $G=\mu_r$, where $r\in\Z^+$. The conjugacy classes
in $W(r,n)=\mu_r\wr S_n$ are in bijection with the set of collections
$(a_{i}(\zeta))_{i\in\Z^+,\zeta\in\mu_r}$
of nonnegative integers such that
\begin{equation} \label{sumeqn}
\sum_{\substack{i\in\Z^+\\\zeta\in\mu_r}} ia_{i}(\zeta) =n.
\end{equation}
In this bijection, $a_{i}(\zeta)$ counts the number of cycles of length $i$
and type $\zeta$. To make this and other statements true for $n=0$, 
we define $W(r,0)$ to be the trivial group.
Following \cite[loc.~cit., \S5]{macdonald}, 
we define the $\C$-algebra
\[ \Lambda(r):=\C[\,p_i(\zeta)\,|\,i\in\Z^+,\zeta\in\mu_r], \]
where $p_i(\zeta)$ are indeterminates.
If $r=1$, we write
$p_i$ instead of $p_{i}(1)$, and identify $\Lambda(1)$ with the complexified
ring of symmetric functions, viewing $p_i$ as the $i$th power sum.
We have an obvious $\N$-grading $\Lambda(r)=\oplus_{n\geq 0}\Lambda(r)_n$,
defined by the rule $\deg(p_i(\zeta))=i$.
If the conjugacy class of $w\in W(r,n)$
corresponds to $(a_{i}(\zeta))$, set
\[ p_w=\prod_{\substack{i\in\Z^+\\\zeta\in\mu_r}} p_i(\zeta)^{a_i(\zeta)}
\ \in \Lambda(r)_n. \]
Let $R(W(r,n))$ denote the space of class functions on $W(r,n)$. If
$f\in R(W(r,n))$, define its \emph{characteristic} by
\[ \ch_{W(r,n)}(f):=\frac{1}{r^n n!}\sum_{w\in W(r,n)}f(w)\,p_w. \]
Then we have an isomorphism $\ch_{W(r,n)}:R(W(r,n))\isomto
\Lambda(r)_n$.
We will also write $\ch_{W(r,n)}(M)$ where $M$ is a
finite-dimensional representation of $W(r,n)$, meaning by this the
characteristic of the character of $M$. (In other words, we will identify
$R(W(r,n))$ with the complexified Grothendieck group of the category of
representations.)
If we combine the maps $\ch_{W(r,n)}$ for all $n\geq 0$,
the resulting isomorphism 
\begin{equation} \label{chareqn}
\ch:\bigoplus_{n\geq 0}R(W(r,n))
\isomto\Lambda(r)
\end{equation}
puts an algebra structure on $\bigoplus_{n\geq 0}R(W(r,n))$,
which is given by induction product (\cite[loc.~cit., (6.3)]{macdonald}).

In order to be able to consider a collection of representations
of $W(r,n)$ for all $n$ simultaneously, we define
the formal power series ring 
\[ \Aa(r):=\C[\![\,p_i(\zeta)\,|\,i\in\Z^+,\zeta\in\mu_r]\!]. \]
Clearly we can extend the isomorphism \eqref{chareqn} to
\begin{equation} \label{chisomeqn}
\ch:\prod_{n\geq 0}R(W(r,n))
\isomto\Aa(r).
\end{equation}
Thus $\Aa(r)$ is a complete $\N$-filtered topological $\C$-algebra.

When dealing with graded representations of $W(r,n)$, we will work in
the ring $\Aa(r)[q]=\Aa(r)\otimes\C[q]$, where
$q$ is a further indeterminate which keeps track of the grading
of the representation.
We give $\Aa(r)[q]$ the $\N$-filtration induced from the
above filtration on $\Aa(r)$ (so $\deg(q)=0$). We will write
$\Aa(r)_+$ and $\Aa(r)[q]_+$ for the augmentation ideals.

For any $P\in\Aa(r)[q]$, we will write $P^\noneq$ for the power series
in $\C[q][\![x]\!]$ obtained by applying the specialization
\begin{equation}\label{specialeqn}
p_{1}(1)\to x, \text{ all other } p_{i}(\zeta)\to 0.
\end{equation}
(Strictly speaking, this rule defines a homomorphism
$\Lambda(r)[q]\to\C[q][x]$, which we then extend to $\Aa(r)[q]$
by continuity. Such an extension will be implicit in similar situations
below.)
Clearly $\ch_{W(r,n)}(f)^\noneq=\frac{1}{r^n n!}f(1)x^n$ for $f\in R(W(r,n))$,
and if $M$ is a representation of $W(r,n)$, $\ch_{W(r,n)}(M)^\noneq=
\frac{1}{r^n n!}(\dim M)x^n$.
So $P^\noneq$ may be thought of as the ``non-equivariant version''
of $P$.

Lastly, we need to introduce plethysm.
Plethysm as defined in \cite[Chapter I, \S8]{macdonald} 
is an associative operation $\circ:\Lambda(1)\times\Lambda(1)
\to\Lambda(1)$, which under the specialization \eqref{specialeqn}
becomes simply the substitution of one element of $\C[x]$ in another.
It is uniquely defined by the following properties:
\begin{enumerate}
\item for all $g\in \Lambda(1)$, the map 
$\Lambda(1)\to \Lambda(1):f\mapsto f\circ g$
is a $\C$-algebra homomorphism;
\item for any $i\in\Z^+$, the map
$\Lambda(1)\to \Lambda(1):g\mapsto p_{i}\circ g$ is a $\C$-algebra
homomorphism;
\item $p_i\circ p_j=p_{ij}$.
\end{enumerate}
Note that $p_1$ is the plethystic identity.
The significance of this operation stems from its interpretation in
terms of polynomial functors, which we will discuss in \S4.

When $r\geq 2$, there is no way to define a similar 
associative operation on $\Lambda(r)$,
but instead there are two operations $\circ:\Lambda(1)\times\Lambda(r)\to
\Lambda(r)$ and $\circ:\Lambda(r)\times\Lambda(1)\to\Lambda(r)$, which
we will also call plethysm. These
have all the associativities one might expect, i.e.\ $f\circ(g\circ h)=
(f\circ g)\circ h$ whenever both sides are defined. The first
operation, which we will not in fact use in the main proof, is defined by
three properties similar to the above, with (3) reading instead 
$p_i\circ p_j(\zeta)
=p_{ij}(\zeta)$. From another viewpoint, this 
is simply the usual kind of plethysm applied to the $r$ tensor factors 
$\C[p_i(\zeta)\,|\,i\in\Z^+]$ of $\Lambda(r)$ independently.
The second operation, which does not seem to have appeared previously
in the literature, is defined by the following properties:
\begin{enumerate}
\item for all $h\in \Lambda(1)$, the map 
$\Lambda(r)\to \Lambda(r):g\mapsto g\circ h$
is a $\C$-algebra homomorphism;
\item for any $i\in\Z^+$, $\zeta\in\mu_r$, the map
$\Lambda(1)\to \Lambda(r):h\mapsto p_{i}(\zeta)\circ h$ is a 
$\C$-algebra homomorphism;
\item $p_i(\zeta)\circ p_j=p_{ij}(\zeta^j)$.
\end{enumerate}
Again, we will give interpretations of these 
operations in terms of polynomial functors in \S4.

For either $r=1$ or $r\geq 2$, the operation
$\circ:\Lambda(r)\times\Lambda(1)\to\Lambda(r)$
can be extended to $\circ:\Aa(r)\times\Aa(1)_+
\to\Aa(r)$. (The assumption of zero constant term on the right is
needed for convergence.) We will need the further
extension of this to an operation 
$\circ:\Aa(r)[q]\times\Aa(1)[q]_+\to\Aa(r)[q]$, which is uniquely defined by:
\begin{enumerate}
\item for all $g\in \Aa(1)[q]_+$, the map 
$\Aa(r)[q]\to \Aa(r)[q]:f\mapsto f\circ g$
is a continuous $\C[q]$-algebra homomorphism;
\item for any $i\in\Z^+$, $\zeta\in\mu_r$, the map
$\Aa(1)[q]_+\to \Aa(r)[q]:g\mapsto p_{i}(\zeta)\circ g$ is a 
continuous $\C$-algebra homomorphism;
\item $p_{i}(\zeta)\circ q=q^i$, $p_i(\zeta)\circ p_j=p_{ij}(\zeta^j)$.
\end{enumerate}
It is easy to see that
for $f\in\Aa(r)[q]$, $g\in\Aa(1)[q]_+$,
\begin{equation} \label{subseqn}
(f\circ g)^\noneq=f^\noneq(g^\noneq),
\end{equation}
where the right-hand side means that $g^\noneq$ is substituted for $x$
in $f^\noneq$.
%%%%%%%%%%%%%%%%%%%%%%%%%%%%%%%%%%%%%%%%%%%%%%%%%%%%%%%%
\section{Statement of Results}
Now we return to the representations of $W(r,n)$ on cohomology discussed
in the Introduction. It is more convenient to state the results in terms
of cohomology
with compact supports. Of course this makes no difference in the case
of the projective variety $\Mbar(r,n)$; and since $\M(r,n)$ is a nonsingular
irreducible variety, Poincar\'e duality gives
\[ H_c^s(\M(r,n),\C)\cong H^{2\dim\M(r,n)-s}(\M(r,n),\C)^* \]
as representations of $W(r,n)$, so the character of the graded representation
$H_c^\bullet(\M(r,n),\C)$ is easily obtained from that of
$H^\bullet(\M(r,n),\C)$ and vice versa.

We will encode the characters of the representations
$H_c^\bullet(\M(r,n),\C)$
and $H^\bullet(\Mbar(r,n),\C)$ for all $n$
in the following elements of $\Aa(r)[q]_+$:
\begin{equation*}
\begin{split}
\CalP(r)&:=
\sum_{n\geq 1}
\sum_{s=\dim\M(r,n)}^{2\dim\M(r,n)}(-1)^s\,\ch_{W(r,n)}(H_c^s(\M(r,n),\C))\,
q^{s-\dim\M(r,n)},\\
\overline{\CalP}(r) &:= \sum_{n\geq 1}
\sum_{s=0}^{\dim\Mbarrn}\ch_{W(r,n)}(H^{2s}(\Mbarrn,\C))\,q^{s}.
\end{split}
\end{equation*}
The non-equivariant versions encode the
Betti numbers:
\begin{equation*}
\begin{split}
\CalP(r)^\noneq&=
\sum_{n\geq 1}\frac{x^n}{r^n n!}
\sum_{s=\dim\M(r,n)}^{2\dim\M(r,n)}
(-1)^s\dim H_c^s(\M(r,n),\C)\,q^{s-\dim\M(r,n)},\\
\overline{\CalP}(r)^\noneq &= \sum_{n\geq 1}\frac{x^n}{r^n n!}
\sum_{s=0}^{\dim\Mbarrn}\dim H^{2s}(\Mbarrn,\C)\,q^{s}.
\end{split}
\end{equation*}
The powers of $q$ used here derive from some easy mixed Hodge theory;
in the terminology of \cite{dimcalehrer},
the terms of these series are the (equivariant or non-equivariant)
\emph{weight polynomials} of the corresponding varieties (with $q=t^2$). 
Explicitly,
since $\M(r,n)$ is either a single point or the complement of 
a nonempty finite set of hyperplanes in a projective space,
it is \emph{minimally pure} in the sense of \cite[Definition 3.1]{dimcalehrer}.
Hence $H_c^s(\M(r,n),\C)$ is nonzero only if 
$\dim\M(r,n)\leq s\leq 2\dim\M(r,n)$,
and is a pure Hodge structure of weight $2s-2\dim\M(r,n)$. By 
\cite[Theorem 5.2]{wonderful},
$\Mbar(r,n)$ has no odd-degree cohomologies, and since it is
nonsingular projective, $H^{2s}(\Mbar(r,n),\C)$ is a pure Hodge structure
of weight $2s$.

In this notation, Lehrer's results are as follows (for $d\in\Z^+$,
$\mu(d)\in\{0,1,-1\}$ is defined in the usual way).
\begin{theorem} \label{1lehrerthm}
In $\Aa(1)[q]$ we have the equation
\[ 1+qp_1+q(q-1)\CalP(1)=\prod_{n\geq 1}(1+p_n)^{R_n/n}, \]
where
\[ R_n:={\sum_{d|n}\mu(d)q^{n/d}}. \]
\end{theorem}
\noindent
This is equivalent to \cite[Theorem 5.5]{lehrerone}; see \S8
for a quick proof due to Getzler.
\begin{theorem} \label{lehrerthm}
In $\Aa(r)[q]$ for $r\geq 2$ we have the equation
\[ 1+(q-1)\CalP(r)=\prod_{\substack{n\geq 1\\\theta\in\mu_r}}
(1+p_{n}(\theta))^{R_{r,n,\theta}/rn}, \]
where
\[ R_{r,n,\theta}:=\sum_{d|n}
|\{\zeta\in\mu_r\,|\,\zeta^d=\theta\}|\,
\mu(d)(q^{n/d}-1). \]
\end{theorem}
\noindent
This can be deduced from \cite[Theorems 5.9 and 6.3]{lehrertwo};
see \S8 for a new proof.
Applying the specialization \eqref{specialeqn}, we find:
\begin{equation} \label{noneqeqn}
\begin{split}
\CalP(1)^\noneq&=\frac{(1+x)^q-1-qx}{q(q-1)},\\
(r\geq 2)\quad 
\CalP(r)^\noneq&=\frac{(1+x)^{\frac{q-1}{r}}-1}{q-1}.
\end{split}
\end{equation}
Of course, these non-equivariant statements can be proved directly:
they are equivalent to special cases of \cite[Theorem 4.8]{orliksolomon}.

Using the familiar $r=1$ kind of plethysm,
the result of Ginzburg and Kapranov
can be expressed as follows.
\begin{theorem} \label{ginzburgthm}
In $\Aa(1)[q]$ we have the equation
\[ \overline{\CalP}(1)= p_1+\CalP(1)\circ\overline{\CalP}(1). \]
Equivalently, $\overline{\CalP}(1)$ is the plethystic inverse of 
$p_1-\CalP(1)$.
\end{theorem}
\noindent
The statement in \cite{ginzburgkapranov} was slightly different.
There is a $\C[q]$-algebra involution $\tau:\Aa(1)[q]\to\Aa(1)[q]$ defined by
\[ \tau(f)=f|_{p_i\to -p_i}, \]
or equivalently by the rule
\[ \tau(\ch_{S_n}(M))=(-1)^{n}\ch_{S_n}(\varepsilon\otimes M), \]
where $M$ is a representation of $S_n$ and $\varepsilon$ is the sign
representation. It is easy to see that $\tau(f\circ g)=\tau(f)\circ(-\tau(g))$,
so Theorem \ref{ginzburgthm} is equivalent to the statement that 
$-\tau(\overline{\CalP}(1))$
is the plethystic inverse of $p_1+\tau(\CalP(1))$.
This follows from \cite[Theorem 3.4.11]{ginzburgkapranov},
bearing in mind \cite[Theorem 3.3.2]{ginzburgkapranov}. 

In \S\!\S 5--7
we will give a simpler proof of Theorem \ref{ginzburgthm}, 
based on hints by Getzler in \cite[\S3]{semiclassical}, in parallel
with the proof of our main result.
Note that Theorem \ref{ginzburgthm} amounts to a recursive 
formula for the terms of $\overline{\CalP}(1)$,
since the degree $n$ terms on the left-hand side
are expressed as polynomials in smaller-degree terms on the right-hand side,
with the known terms of $\CalP(1)$ as coefficients. (Here we need
to observe that the degree $1$ term of $\CalP(1)$ vanishes,
since $\M(1,1)$ is empty.)
Using \eqref{subseqn} and \eqref{noneqeqn}, we get the
non-equivariant version:
\begin{equation}
\overline{\CalP}(1)^\noneq = x+\frac{(1+\overline{\CalP}(1)^\noneq)^q-
1-q\overline{\CalP}(1)^\noneq}{q(q-1)}.
\end{equation}
This is equivalent to Keel's recursion \cite[p.~550, (4)]{keel},
and to \cite[(0.7)]{manin}.

Using the generalized kind of plethysm, we
can state the main result of this paper, which will be proved in \S\!\S 5--7.
\begin{theorem} \label{mainthm}
In $\Aa(r)[q]$ for $r\geq 2$ we have the equation:
\[ \overline{\CalP}(r)=(1+\overline{\CalP}(r))
(\CalP(r)\circ\overline{\CalP}(1)). \]
Equivalently, $1+\overline{\CalP}(r)$ is the multiplicative
inverse of $1-\CalP(r)\circ\overline{\CalP}(1)$.
\end{theorem}
\noindent
Of course, this amounts to a formula for the degree $n$ term of
$\overline{\CalP}(r)$ in terms of the degree $\leq n$ terms of
$\CalP(r)\circ\overline{\CalP}(1)$, which can be calculated using
Theorems \ref{lehrerthm} and \ref{ginzburgthm}. Using \eqref{subseqn}
and \eqref{noneqeqn},
we get the non-equivariant version:
\begin{corollary} $(r\geq 2)$
In $\C[q][\![x]\!]$, $1+\overline{\CalP}(r)^\noneq$ is the multiplicative
inverse of $1-\frac{1}{q-1}((1+\overline{\CalP}(1)^\noneq)^{\frac{q-1}{r}}-1)$.
\end{corollary}
\noindent
The special case $r=2$ of this Corollary was obtained by Yuzvinsky
by constructing explicit bases for the cohomology groups.
See \cite[Theorem 5.1(i)]{yuzvinsky}; his $t$ is our $x$, and
his $\alpha$ is our $(1+\overline{\CalP}(1)^\noneq)^{\frac{1}{2}}$.
%%%%%%%%%%%%%%%%%%%%%%%%%%%%%%%%%%%%%%%%%%%%%%%%%%%%%%%%
\section{An Example}
Let us calculate $\overline{\CalP}(2)$ up to degree $3$, using Theorem
\ref{mainthm}. To make the formulas more legible, 
we will write $x_i$ instead of 
$p_{i}(1)$ and $y_i$ instead of $p_i(-1)$. Thus $\Aa(2)[q]=
\C[q][\![x_1,y_1,x_2,y_2,\cdots]\!]$, and plethysm rule (3) becomes
\begin{equation*}
\begin{split}
&x_i\circ q = y_i\circ q=q^i,\\
&x_i\circ p_j=x_{ij},\\
&y_i\circ p_j=\left\{\begin{array}{cl}
y_{ij},&\text{ if $j$ is odd}\\
x_{ij},&\text{ if $j$ is even.}
\end{array}
\right.
\end{split}
\end{equation*}
Now $\CalP(2)$ is given by Theorem \ref{lehrerthm}.
Since we are concerned with degree $\leq 3$, we only need those
factors of the infinite product where $n\leq 3$. These factors are:
\[ (1+x_1)^{\frac{q-1}{2}}(1+y_1)^{\frac{q-1}{2}}
(1+x_2)^{\frac{q^2-2q+1}{4}}(1+y_2)^{\frac{q^2-1}{4}}
(1+x_3)^{\frac{q^3-q}{6}}(1+y_3)^{\frac{q^3-q}{6}}. \]
So up to degree $3$,
\begin{equation} \label{p2eqn}
\begin{split}
\CalP(2)&=\frac{1}{2}(x_1+y_1)\\
&+\frac{q-3}{8}(x_1^2+y_1^2)
+\frac{q-1}{4}x_1y_1+\frac{q-1}{4}x_2+\frac{q+1}{4}y_2\\
&+\frac{q^2-8q+15}{48}(x_1^3+y_1^3)+\frac{q^2-4q+3}{16}(x_1^2y_1+x_1y_1^2)\\
&\qquad+\frac{q^2-1}{8}(x_1y_2+y_1y_2)+\frac{q^2-2q+1}{8}(x_1x_2+x_2y_1)\\
&\qquad\qquad+\frac{q^2+q}{6}(x_3+y_3)+\cdots
\end{split}
\end{equation}
The denominators of the coefficients are the orders of the
centralizers of the corresponding elements of $W(2,n)$, so the numerators
represent the actual traces on the graded cohomology. Since the centre
$\{1,-1\}$ of $W(2,n)$ acts trivially on $H_c^\bullet(\M(2,n))$, this series
is stable under the involution corresponding to multiplication by $-1$,
which is just
\begin{equation*}
\begin{split}
(i \text{ odd})\quad  &x_i\mapsto y_i, y_i\mapsto x_i\\
(i \text{ even})\quad &x_i\mapsto x_i, y_i\mapsto y_i.
\end{split}
\end{equation*}
The terms have been grouped to reflect this. Now up to degree $3$,
\begin{equation} \label{pbar1eqn}
\overline{\CalP}(1)=p_1+\frac{1}{2}(p_1^2+p_2)
+\frac{q+1}{6}p_1^3+\frac{q+1}{2}p_1p_2+\frac{q+1}{3}p_3+\cdots
\end{equation}
(Thus far the series $\overline{\CalP}(1)$ involves only
trivial representations of the symmetric groups; to calculate further terms,
one could use the recursion given by Theorem \ref{ginzburgthm}.)
Applying the plethysm rules to \eqref{p2eqn} and \eqref{pbar1eqn} 
tells us that
\begin{equation} \label{inveqn}
\begin{split}
\CalP(2)\circ\overline{\CalP}(1)&=\frac{1}{2}(x_1+y_1)\\
&+\frac{q-1}{8}(x_1^2+y_1^2)+\frac{q-1}{4}x_1y_1+\frac{q+1}{4}x_2
+\frac{q+1}{4}y_2\\
&+\frac{q^2+2q+1}{48}(x_1^3+y_1^3)+\frac{q^2-2q+1}{16}(x_1^2y_1+x_1y_1^2)\\
&\qquad+\frac{q^2-1}{8}(x_1y_2+y_1y_2)+\frac{q^2+2q-1}{8}(x_1x_2+x_2y_1)\\
&\qquad\qquad+\frac{q^2+2q+1}{6}(x_3+y_3)+\cdots
\end{split}
\end{equation}
By Theorem \ref{mainthm}, $1+\overline{\CalP}(2)$
is the inverse of $1-\CalP(2)\circ\overline{\CalP}(1)$
in $\Aa(2)[q]$. Hence we find
\begin{equation}
\begin{split}
\overline{\CalP}(2)&=\frac{1}{2}(x_1+y_1)\\
&+\frac{q+1}{8}(x_1^2+y_1^2)+\frac{q+1}{4}x_1y_1+\frac{q+1}{4}x_2
+\frac{q+1}{4}y_2\\
&+\frac{q^2+8q+1}{48}(x_1^3+y_1^3)+\frac{q^2+4q+1}{16}(x_1^2y_1+x_1y_1^2)\\
&\qquad+\frac{q^2+2q+1}{8}(x_1y_2+y_1y_2)+\frac{q^2+4q+1}{8}(x_1x_2+x_2y_1)\\
&\qquad\qquad+\frac{q^2+2q+1}{6}(x_3+y_3)+\cdots
\end{split}
\end{equation}
The only nontrivial representation involved here is $H^2(\overline{\M(2,3)})$.
This calculation shows that the traces of various elements of $W(2,3)$
on this cohomology are $8$, $4$, and $2$. Since $\{\pm 1\}$ acts trivially,
we may regard this as a representation of $W(2,3)/\{\pm 1\}\cong S_4$.
Consulting the character table of $S_4$, we find that its isomorphism type is
$3V_{(4)}\oplus V_{(31)}\oplus V_{(2^2)}$, where the
irreducible representations of $S_4$ are parametrized by
partitions of $4$ as usual.
%%%%%%%%%%%%%%%%%%%%%%%%%%%%%%%%%%%%%%%%%%%%%%%%%%%%%%%%%%%
\section{Species, $\B_r$-modules, and Polynomial Functors}
Before proceeding to the proof of Theorems \ref{ginzburgthm} and
\ref{mainthm}, we need to build more of the combinatorial framework.
As above, let $r\in\Z^+$.
We introduce the category $\B_r$,
whose objects are finite sets with a free action of $\mu_r$, 
and whose morphisms are bijections respecting this action.
(Clearly $\B_1$ is just the category of finite sets and bijections.)
For all $n\in\N$, write $[n]$ for $\{1,\cdots,n\}$ ($[0]$ is the empty set).
The set $[n]_r:=\mu_r\times [n]$,
with the obvious $\mu_r$-action by left multiplication on the first factor,
is an object of $\B_r$, and clearly all objects of $\B_r$ are isomorphic
to $[n]_r$ for unique $n\in\N$.
The group of automorphisms of the object
$[n]_r$ in the category $\B_r$ is canonically isomorphic to the wreath
product $W(r,n)$. (This is true when $n=0$,
since we have defined $W(r,0)$ to be the trivial group.)

In enumerative combinatorics (for instance, see \cite{joyalsln} or \cite{bll}) 
a \emph{species} is a functor $A:\B_1\to\B_1$. More generally, we will
define an \emph{$r$-species} to be a functor $A:\B_r\to\B_1$.
Note that $A$ is determined up to (non-unique) isomorphism by the
collection of sets $A([n]_r)$ for $n\geq 0$, each equipped with the
action of $W(r,n)$ induced by $A$. We
define the \emph{cycle index series} $Z_A\in\Aa(r)$ by
\[ Z_A:=\sum_{n\geq 0}\frac{1}{r^n n!}\sum_{w\in W(r,n)}
|A([n]_r)^{A(w)}|\, p_w. \]
This generalizes Joyal's definition in the case $r=1$. 

For example, we have a ``trivial'' $r$-species $E(r)$ defined by
$E(r)(I)=\{I\}$ for all $I$ in $\B_r$. (The definition of $E(r)$ on 
morphisms is forced.) For any $n\in\N$, we define an $r$-species $E(r)_n$ to be
the degree-$n$ sub-$r$-species of $E(r)$. That is,
\[ E(r)_n(I)=\left\{\begin{array}{cl}
\{I\},&\text{ if $I\cong [n]_r$ in $\B_r$}\\
\emptyset,&\text{ otherwise.}
\end{array}\right. \]
Similarly we define $E(r)_+$, which takes the value
$\{I\}$ when $I\neq\emptyset$ and $\emptyset$ otherwise,
and $E(r)_{\geq 2}$, which takes the value $\{I\}$ when $|I|\geq 2r$
and $\emptyset$ otherwise.
We have
\begin{equation} \label{expcheqn}
Z_{E(r)}=\sum_{n\geq 0}\frac{1}{r^n n!}{\sum_{w\in W(r,n)}}p_w
=\exp(\sum_{i\geq 1}\frac{1}{ri}\sum_{\zeta\in\mu_r} p_i(\zeta)),
\end{equation}
where the second equality follows from the fact that an element of
$W(r,n)$ whose class corresponds to $(a_i(\zeta))$ has centralizer of order
\begin{equation} \label{centordereqn}
\prod_{\substack{i\geq 1\\\zeta\in\mu_r}} a_i(\zeta)! (ri)^{a_i(\zeta)}. 
\end{equation}
Clearly $Z_{E(r)_+}$ and $Z_{E(r)_{\geq 2}}$ are obtained from
$Z_{E(r)}$ by subtracting the degree $0$ and the degree $\leq 1$ terms
respectively.

There is a linear analogue of the concept
of $r$-species. A \emph{$\B_r$-module} is a functor
$\B_r\to\Cmod$, where $\Cmod$
is the category of finite-dimensional complex vector spaces.
(A $\B_1$-module
is what is elsewhere called a \emph{tensor species} or \emph{$\Ss$-module}.)
A $\B_r$-module $U$ is determined up to isomorphism by
the vector spaces $(U([n]_r))$ for $n\in\N$, each 
equipped with the representation of $W(r,n)$ induced by $U$.
Recalling the isomorphism \eqref{chisomeqn}, it is natural to define
the \emph{characteristic} $\ch(U)\in\Aa(r)$ by
\[ \ch(U):=\sum_{n\geq 0} \ch_{W(r,n)}(U([n]_r)). \]

Let $H^0:\B_1\to\Cmod$ be the functor defined by
$H^0 I:=\C^I$, with the obvious definition on morphisms.
By composing with $H^0$, any $r$-species $A:\B_r\to\B_1$ gives rise
to a $\B_r$-module $H^0 A:\B_r\to\Cmod$, which may be thought of
as the ``linearization'' of $A$. Moreover, it is clear that $\ch(H^0 A)=Z_A$,
so many facts about cycle index series of $r$-species
can be deduced from more general results about characteristics of
$\B_r$-modules. We will take this approach, instead of working with
$r$-species directly.

The main point of Joyal's theory is that various natural operations
on species (or tensor species) correspond to natural operations on
their cycle index series (or characteristics). This remains true
for general $r$.
For instance, we define the sum and product of two $r$-species 
$A$ and $B$ as follows:
\begin{equation}
\begin{split}
(A+B)(I)&=A(I)\amalg B(I),\\
(A\cdot B)(I)&=\coprod_{(I_1,I_2)\in\Decomp(I)^{\mu_r}}
A(I_1)\times B(I_2).
\end{split}
\end{equation}
Here $\Decomp(I)^{\mu_r}$ is the set of all decompositions of $I$ into
the disjoint union of two (ordered) $\mu_r$-stable subsets.
The definitions on morphisms of $\B_r$ are obvious.
Analogously, the sum and product of two $\B_r$-modules $U$ and $V$ are defined
as follows:
\begin{equation} \label{sumproddefeqn}
\begin{split}
(U+V)(I)&=U(I)\oplus V(I),\\
(U\cdot V)(I)&=\bigoplus_{(I_1,I_2)\in\Decomp(I)^{\mu_r}}
U(I_1)\otimes V(I_2).
\end{split}
\end{equation}
Note that this is just the induction product on
$\prod_{n\geq 0}R(W(r,n))$, in the sense that
\[ (U\cdot V)([n]_r)\cong\bigoplus_{n_1+n_2=n}
\Ind_{W(r,n_1)\times W(r,n_2)}^{W(r,n)}(U([n_1]_r)\otimes V([n_2]_r)) \]
as representations of $W(r,n)$.
\begin{proposition} \label{sumprodprop}
If $U,V$ are $\B_r$-modules, then
\begin{equation*}
\ch(U+V)=\ch(U)+\ch(V),\ \ch(U\cdot V)=\ch(U)\ch(V).
\end{equation*}
\end{proposition}
\begin{proof}
The first statement is obvious, and the second follows from
\cite[Chapter I, Appendix B, (6.3)]{macdonald}.
\end{proof}
\begin{corollary} \label{specsumprodcor}
If $A,B$ are $r$-species, then
\[ Z_{A+B}=Z_A+Z_B,\ Z_{A\cdot B}=Z_A Z_B. \]
\end{corollary}
\begin{proof}
Since $H^0(A+B)\cong H^0 A + H^0 B$ and similarly for multiplication,
this follows directly from Proposition \ref{sumprodprop}.
\end{proof}
\noindent
In the $r=1$ case, these results are well known (see \cite[\S1.3]{bll}
or \cite{joyalsln}).

More significantly, we have operations of \emph{substitution} (or
\emph{partitional composition}) of
$r$-species and $\B_r$-modules which correspond to the operations of plethysm
introduced in \S1. In the $r=1$ case, Joyal defined the substitution
$A\circ B$ of two species $A$ and $B$ to be the species whose value on
an object $I$ of $\B_1$ is
\begin{equation}
(A\circ B)(I)=\coprod_{\pi\in\Par(I)}
\left(A(\pi)\times\prod_{J\in\pi}B(J)\right),
\end{equation}
where $\Par(I)$ is the set of partitions of $I$, i.e.\
sets of non-empty disjoint subsets whose union equals $I$. The definition
on morphisms is natural. (Note that $B(\emptyset)$ is irrelevant to this
definition, and it is common to assume that it is empty.) There are two
generalizations of this to the context of $r$-species.

First, if $A$ is an ordinary species and
$B$ is an $r$-species, one can define an $r$-species $A\circ B$ whose
value on an object $I$ of $\B_r$ is
\begin{equation} 
(A\circ B)(I)=\coprod_{\pi\in\Par(I)_{\mu_r}}
\left(A(\pi)\times\prod_{J\in\pi}B(J)\right).
\end{equation}
Here $\Par(I)_{\mu_r}$ is the set of partitions $\pi\in\Par(I)$ for
which $\mu_r$ stabilizes each part. Analogously, if $U$ is a $\B_1$-module
and $V$ is a $\B_r$-module, we define a $\B_r$-module $U\circ V$ by
\begin{equation}
(U\circ V)(I)=\bigoplus_{\pi\in\Par(I)_{\mu_r}}
\left(U(\pi)\otimes\bigotimes_{J\in\pi}V(J)\right).
\end{equation}
As we will see below, this corresponds to the less interesting 
kind of plethysm.

Second, if $A$ is an $r$-species and
$B$ is an ordinary species, one can define an $r$-species $A\circ B$ whose
value on an object $I$ of $\B_r$ is
\begin{equation}
(A\circ B)(I)=\coprod_{\pi\in\Par(I)_{\B_r}}
\left(A(\pi)\times\prod_{\CalO\in\mu_r\setminus\pi}B(\CalO)\right).
\end{equation}
Here $\Par(I)_{\B_r}$ is the set of partitions which are preserved by the 
$\mu_r$-action on $I$, and such that the $\mu_r$-action
on $\pi$ (i.e.\ the set of parts of the partition) is free, so that
$\pi$ is an object of $\B_r$.
For every $\mu_r$-orbit $\CalO\subseteq\pi$, $B(\CalO)$ is the limit
of $B(J)$ for $J\in\CalO$; in other words,
\[ B(\CalO):=\{(b_J)\in\prod_{J\in\CalO} B(J)\,|\,
b_{\zeta.J}=B(\zeta)(b_J),\ \forall J\in\CalO, \zeta\in\mu_r\}, \]
where $B(\zeta):B(J)\to B(\zeta.J)$ is the bijection induced
by the bijection $\zeta:J\to \zeta.J$. 
Analogously, if $U$ is a $\B_r$-module and $V$ is a $\B_1$-module, we define
a $\B_r$-module $U\circ V$ by
\begin{equation} 
(U\circ V)(I)=\bigoplus_{\pi\in\Par(I)_{\B_r}}
\left(U(\pi)\otimes\bigotimes_{\CalO\in\mu_r\setminus\pi}V(\CalO)\right).
\end{equation}
As we will see below, this corresponds to the more interesting kind
of plethysm.

We say that a $\B_r$-module $U$ is \emph{bounded} if $U(I)=\{0\}$ for 
$|I|$ sufficiently large, or equivalently if $\ch(U)\in\Lambda(r)$.
For any $\B_r$-module $U$, we can define bounded ``truncations''
$U_{\leq N}$ by
\[ U_{\leq N}(I)=\left\{\begin{array}{cl}
U(I),&\text{ if $|I|\leq rN$,}\\
\{0\},&\text{ if $|I|> rN$.}
\end{array}\right. \]
Clearly $\ch(U)=\lim_{N\to\infty}\ch(U_{\leq N})$, which means that
some properties of characteristics need only be proved for bounded
$\B_r$-modules.

Let $\C\mu_r$-$\mathbf{mod}$ be the category of
finite-dimensional representations of $\mu_r$.
Any bounded $\B_r$-module $U$ gives rise to a 
functor $F_U:\C\mu_r$-$\mathbf{mod}\to
\C$-$\mathbf{mod}$, defined on objects by
\[ F_U(M):=\bigoplus_{n\geq 0}(U([n]_r)\otimes M^{\otimes n})^{W(r,n)}. \]
Here the action of $W(r,n)$ on $M^{\otimes n}$ is given by
the permutation of the tensor factors composed with the actions
of $\mu_r$ on each. The functor $F_U$ is \emph{polynomial} in the sense
that for all $M,N\in\C\mu_r$-$\mathbf{mod}$, the map
$F_U:\Hom_{\C\mu_r}(M,N)\to
\Hom(F_U(M),F_U(N))$ is a polynomial function. 
It is a consequence of \cite[Theorem 6.4]{polyfunctors} that all polynomial
functors $\C\mu_r$-$\mathbf{mod}\to
\C$-$\mathbf{mod}$ are isomorphic to one of the form $F_U$,
where $U$ is a bounded $\B_r$-module.

Without the boundedness assumption, one can define an \emph{analytic}
functor in the same way (relaxing the finite-dimensionality condition
on the vector spaces involved); see \cite[Chap.~4]{joyalsln} for the 
$r=1$ case. The reason for imposing the boundedness
assumption is that we then have a simple
relationship between $\ch(U)$ and $F_U$:
\begin{proposition} \label{traceprop}
If $U$ is a bounded $\B_r$-module, $M\in\C\mu_r$-$\mathbf{mod}$,
and $\varphi\in\End(M)$ commutes with the $\mu_r$-action, then
\[ \ch(U)|_{p_i(\zeta)\to\tr(\varphi^i\zeta,M)}=\tr(F_U(\varphi),F_U(M)). \]
\end{proposition}
\begin{proof}
It is well known that if $N$ is a finite-dimensional representation
of a finite group $\Gamma$ and 
$\psi\in\End(N)$ commutes with the $\Gamma$-action, then
\begin{equation} 
\tr(\psi,N^{\Gamma})=\frac{1}{|\Gamma|}\sum_{\gamma\in\Gamma}
\tr(\gamma\psi,N).
\end{equation}
Applying this with $\Gamma=W(r,n)$, $N=U([n]_r)\otimes M^{\otimes n}$
and $\psi=\mathrm{id}\otimes\varphi^{\otimes n}$, we get that
\[ \tr(F_U(\varphi),F_U(M))=\sum_{n\geq 0}\frac{1}{|W(r,n)|}
\sum_{w\in W(r,n)}\negthickspace\negthickspace
\tr(w,U([n]_r))\tr(w\varphi^{\otimes n},M^{\otimes n}). \]
So we need only show that for $w\in W(r,n)$,
\[ p_w|_{p_i(\zeta)\to\tr(\varphi^i\zeta,M)}=
\tr(w\varphi^{\otimes n},M^{\otimes n}). \]
This clearly reduces to the case when $w$ consists of a single
cycle of length $n$ and type $\zeta$, in which case it says that
\begin{equation} \label{fundtraceeqn} 
\tr(\varphi^n\zeta,M)=\tr(w\varphi^{\otimes n},M^{\otimes n}).
\end{equation}
This is easy to prove: one way is to observe that it suffices to
prove this for $\varphi$ a regular semisimple element of $\End(M)$,
and then split $M$ into $\varphi$-eigenspaces, reducing
to the trivial case where $\dim M=1$.
\end{proof}
\noindent
In the case $r=1$, this Proposition is equivalent to remarks in
\cite[Chapter I, Appendix A, \S7]{macdonald}; the specialization
sends the ``abstract'' $i$th power sum $p_i$ to the actual sum of the
$i$th powers of the eigenvalues of $\varphi$. For general $r$, it 
is easy to see that
$f\in\Aa(r)$ is uniquely determined by the specializations
$f|_{p_i(\zeta)\to\tr(\varphi^i\zeta,M)}$ for all $M$ and $\varphi$,
so Proposition \ref{traceprop} gives a way of determining $\ch(U)$ from $F_U$
when $U$ is bounded.

Now we must reinterpret our operations on $\B_r$-modules in terms of
polynomial functors. Addition and multiplication correspond to the
obvious ``pointwise'' operations:
\begin{proposition} \label{functsumprodprop}
If $U,V$ are bounded $\B_r$-modules, then
\[ F_{U+V}\cong F_U\oplus F_V,\ F_{U\cdot V}\cong F_U\otimes F_V. \]
\end{proposition}
\begin{proof}
The statement about addition is obvious. The assertion that
$F_{U\cdot V}\cong F_U\otimes F_V$ boils down to the fact that
for any $M\in\C\mu_r$-$\mathbf{mod}$,
\begin{equation*}
\begin{split}
&(U([n_1]_r)\otimes M^{\otimes n_1})^{W(r,n_1)}
\otimes(V([n_2]_r)\otimes M^{\otimes n_2})^{W(r,n_2)}\\
&\cong\left(\Ind_{W(r,n_1)\times W(r,n_2)}^{W(r,n_1+n_2)}
(U([n_1]_r)\otimes V([n_2]_r))\otimes M^{\otimes(n_1+n_2)}\right)
^{W(r,n_1+n_2)},
\end{split}
\end{equation*}
and that this isomorphism is natural in $M$.
This follows from the Frobenius reciprocity rule
\begin{equation} \label{frobeniuseqn}
(M\otimes\mathrm{Res}_H^G N)^H\cong (\Ind_H^G M\otimes N)^G,
\end{equation}
applied with $H=W(r,n_1)\times W(r,n_2)$, $G=W(r,n_1+n_2)$.
\end{proof}
\noindent
Note that we could have deduced Proposition \ref{sumprodprop} from this,
using Proposition \ref{traceprop}.

The upshot of \cite[Chapter I, Appendix A, \S6]{macdonald} is that
substitution of bounded $\B_1$-modules corresponds to composition of
the associated polynomial functors; in \cite[loc.~cit., \S7]{macdonald}
this is used to prove that substitution corresponds to plethysm of 
characteristics. (For a more detailed proof, take the $r=1$ case of Theorems
\ref{functisomthm} and \ref{plethysmthm} below.) This remains true in the case
$r\geq 2$, except that of course it makes no sense to compose two 
functors $\C\mu_r$-$\mathbf{mod}\to
\C$-$\mathbf{mod}$. The two possible types of composition
``explain'' the two generalizations of substitution, and the
two generalizations of plethysm.

Firstly, if $F:\C$-$\mathbf{mod}\to
\C$-$\mathbf{mod}$ and $G:\C\mu_r$-$\mathbf{mod}\to\C$-$\mathbf{mod}$ are 
polynomial functors, then so is $F\circ G$.
\begin{theorem} \label{otherfunctisomthm}
If $U$ is a bounded $\B_1$-module, and $V$ is a bounded $\B_r$-module
with $V(\emptyset)=\{0\}$, then $F_{U\circ V}\cong F_U\circ F_V$.
\end{theorem}
\begin{proof}
This is implicit in \cite[\S9]{polyfunctors}, and can be proved by
an argument very similar to that for Theorem \ref{functisomthm} below.
\end{proof}
\noindent
We can now prove that this kind of substitution corresponds to the first
generalization of plethysm (in the $r=1$ case, this is a result in
\cite[Section 4.4]{joyalsln}).
\begin{theorem}
If $U$ is a $\B_1$-module, and $V$ is a $\B_r$-module with 
$V(\emptyset)=\{0\}$, then $\ch(U\circ V)=\ch(U)\circ\ch(V)$.
\end{theorem}
\begin{proof}
It is easy to see that
$\ch(U\circ V)$ and $\ch(U)\circ\ch(V)$ are the limits as $N\to\infty$
of $\ch(U_{\leq N}\circ V_{\leq N})$ and
$\ch(U_{\leq N})\circ\ch(V_{\leq N})$ respectively. Hence 
it suffices to prove the Theorem
when $U$ and $V$ are bounded. Then, as observed above, it suffices to 
show that $\ch(U\circ V)$ and $\ch(U)\circ\ch(V)$ become equal
under every specialization of the form
$p_i(\zeta)\to\tr(\varphi^i\zeta,M)$, where $M\in\C\mu_r$-$\mathbf{mod}$
and $\varphi\in\End(M)$ commutes with the $\mu_r$-action.
Now using Proposition \ref{traceprop} and Theorem \ref{otherfunctisomthm},
\begin{equation*}
\begin{split}
\ch(U\circ V)|_{p_i(\zeta)\to\tr(\varphi^i\zeta,M)}
&=\tr(F_{U\circ V}(\varphi),F_{U\circ V}(M))\\
&=\tr(F_U(F_V(\varphi)),F_U(F_V(M)))\\
&=\ch(U)|_{p_i\to\tr(F_V(\varphi)^i,F_V(M))}\\
&=\ch(U)|_{p_i\to\tr(F_V(\varphi^i),F_V(M))}\\
&=\ch(U)|_{p_i\to\ch(V)|_{p_j(\zeta)\to\tr(\varphi^{ij}\zeta,M)}}.
\end{split}
\end{equation*}
That this equals $(\ch(U)\circ\ch(V))|_{p_i(\zeta)\to\tr(\varphi^i\zeta,M)}$
is obvious from the definition of this kind of plethysm.
\end{proof}
\begin{corollary}
If $A$ is a species and $B$ is an $r$-species with $B(\emptyset)=\emptyset$,
then $Z_{A\circ B}=Z_A\circ Z_B$.
\end{corollary}
\noindent
We will not actually use these results about the first generalization 
of plethysm in
the remainder of the paper.

Secondly, and more importantly, if $F:\C\mu_r$-$\mathbf{mod}\to
\C$-$\mathbf{mod}$ and $G:\C$-$\mathbf{mod}\to\C$-$\mathbf{mod}$ are 
polynomial functors, then
$G$ induces a functor $G^{(r)}:\C\mu_r$-$\mathbf{mod}\to
\C\mu_r$-$\mathbf{mod}$, and we can compose to get a polynomial functor
$F\circ G^{(r)}:\C\mu_r$-$\mathbf{mod}\to\C$-$\mathbf{mod}$.
\begin{theorem} \label{functisomthm}
If $U$ is a bounded $\B_r$-module, and $V$ is a bounded $\B_1$-module
with $V(\emptyset)=\{0\}$, then
$F_{U\circ V}\cong F_U\circ F_V^{(r)}$.
\end{theorem}
\begin{proof}
Let $M\in\C\mu_r$-$\mathbf{mod}$. We will give a chain
of isomorphisms showing that $F_{U\circ V}(M)\cong F_U(F_V^{(r)}(M))$,
and leave it to the reader to check that the isomorphisms are natural
and hence constitute an isomorphism of functors. By definition,
\[ F_{U\circ V}(M)=\bigoplus_{n\geq 0}\left(
\bigoplus_{\pi\in\Par([n]_r)_{\B_r}}
(U(\pi)\otimes\bigotimes_{\CalO\in\mu_r\setminus\pi}V(\CalO))
\otimes M^{\otimes n}\right)^{W(r,n)}. \]
Now for any $\pi\in\Par([n]_r)_{\B_r}$,
the parts of $\pi$ occur in $\mu_r$-orbits of size $r$. The size
of the part is constant on each orbit; these sizes, arranged in descending 
order, give a partition $n_1\geq n_2\geq \cdots\geq n_m$ of $n$, where
$n_i\in\Z^+$. This partition $(n_i)$ is the \emph{type} of $\pi$.
Fixing $n_1\geq n_2\geq \cdots\geq n_m$, the set of all 
$\pi\in\Par([n]_r)_{\B_r}$ with $\mathrm{type}(\pi)=(n_i)$ form
a single $W(r,n)$-orbit. Let $W(r,(n_i))$ be the stabilizer of the
``standard'' partition $\pi_{n_1,n_2,\cdots,n_m}\in\Par([n]_r)_{\B_r}$,
whose parts are
\[ \{(\zeta,j)\,|\,n_1+\cdots+n_{k-1}+1\leq j\leq n_1+\cdots+n_{k}\},
\ \zeta\in\mu_r, 1\leq k\leq m. \]
Then we have an isomorphism of representations of $W(r,n)$:
\[ \bigoplus_{\substack{\pi\in\Par([n]_r)_{\B_r}\\\mathrm{type}(\pi)=(n_i)}}
\negthickspace\negthickspace\negthickspace\negthickspace
(U(\pi)\otimes\bigotimes_{\CalO\in\mu_r\setminus\pi}V(\CalO))
\cong\Ind_{W(r,(n_i))}^{W(r,n)}(U([m]_r)\otimes V[n_1]\otimes\cdots\otimes
V[n_m]). \]
Here we have identified the set of parts of $\pi_{n_1,n_2,\cdots,n_m}$
with $[m]_r$ in the obvious way, and identified $V(\CalO)$ with $V[n_i]$
if $\CalO$ is the $\mu_r$-orbit on these parts corresponding to the
$\mu_r$-orbit of $(1,i)$. Applying \eqref{frobeniuseqn} with
$H=W(r,(n_i))$, $G=W(r,n)$, we can transform our expression
for $F_{U\circ V}(M)$ into
\[ \bigoplus_{\substack{m\geq 0\\n_1\geq\cdots\geq n_m\geq 1}}
\negthickspace\negthickspace\negthickspace\negthickspace
\left(U([m]_r)\otimes V[n_1]\otimes\cdots\otimes
V[n_m]\otimes M^{\otimes (n_1+\cdots+n_m)}\right)^{W(r,(n_i))}. \]
Now $W(r,(n_i))$ has a normal subgroup, which we will identify
with $S_{n_1}\times\cdots\times S_{n_m}$, consisting of all $w\in W(r,n)$
which stabilize each individual part of $\pi_{n_1,n_2,\cdots,n_m}$.
This subgroup acts trivially on $U([m]_r)$. Let $\overline{W}(r,(n_i))$
be the quotient group $W(r,(n_i))/(S_{n_1}\times\cdots\times S_{n_m})$.
Using the obvious fact that $(-)^{W(r,(n_i))}=
((-)^{S_{n_1}\times\cdots\times S_{n_m}})^{\overline{W}(r,(n_i))}$,
we can rewrite the previous expression as
\[ \bigoplus_{\substack{m\geq 0\\n_1\geq\cdots\geq n_m\geq 1}}
\negthickspace\negthickspace\negthickspace\negthickspace
\left(U([m]_r)\otimes (V[n_1]\otimes M^{\otimes n_1})^{S_{n_1}}
\otimes\cdots\otimes
(V[n_m]\otimes M^{\otimes n_m})^{S_{n_m}}
\right)^{\overline{W}(r,(n_i))}. \]
Now $\overline{W}(r,(n_i))$ can be identified with the subgroup
of $W(r,m)$ consisting of all $w\in W(r,m)$ whose image $\tilde{w}\in S_m$
satisfies $n_{\tilde{w}(i)}=n_i$ for all $i\in[m]$. The action of
$\overline{W}(r,(n_i))$ on $U([m]_r)$ is just the restriction of the
action of $W(r,m)$. Since
\begin{equation*}
\begin{split} 
\Ind_{\overline{W}(r,(n_i))}^{W(r,m)}&\left(
\bigoplus_{n_1\geq\cdots\geq n_m\geq 1}\negthickspace\negthickspace
(V[n_1]\otimes M^{\otimes n_1})^{S_{n_1}}\otimes\cdots\otimes
(V[n_m]\otimes M^{\otimes n_m})^{S_{n_m}}\right)\\
&\cong\bigoplus_{n_1,\cdots,n_m\geq 1}\negthickspace
(V[n_1]\otimes M^{\otimes n_1})^{S_{n_1}}\otimes\cdots\otimes
(V[n_m]\otimes M^{\otimes n_m})^{S_{n_m}}\\
&\cong(\bigoplus_{n\geq 1}(V[n]\otimes M^{\otimes n})^{S_n})^{\otimes m},
\end{split}
\end{equation*}
we can apply \eqref{frobeniuseqn} again to transform our expression
for $F_{U\circ V}(M)$ into
\[ \bigoplus_{m\geq 0}\left(U([m]_r)\otimes
(\bigoplus_{n\geq 1}(V[n]\otimes M^{\otimes n})^{S_n})^{\otimes m}
\right)^{W(r,m)}, \]
which is $F_U(F_V^{(r)}(M))$ by definition.
\end{proof}

Now we come to the main result of this section
(another generalization of Joyal's result in
\cite[Section 4.4]{joyalsln}), which shows that the second kind of substitution
corresponds to the more interesting kind of plethysm.
\begin{theorem} \label{plethysmthm}
If $U$ is a $\B_r$-module, and $V$ is a $\B_1$-module with
$V(\emptyset)=\{0\}$, then $\ch(U\circ V)=\ch(U)\circ\ch(V)$.
\end{theorem}
\begin{proof}
Since $\ch(U\circ V)$ and $\ch(U)\circ\ch(V)$ are the limits as $N\to\infty$
of $\ch(U_{\leq N}\circ V_{\leq N})$ and
$\ch(U_{\leq N})\circ\ch(V_{\leq N})$ respectively, it suffices to prove this 
when $U$ and $V$ are bounded. Then, as before, it suffices to 
show that $\ch(U\circ V)$ and $\ch(U)\circ\ch(V)$ become equal
under every specialization of the form
$p_i(\zeta)\to\tr(\varphi^i\zeta,M)$, where $M\in\C\mu_r$-$\mathbf{mod}$
and $\varphi\in\End(M)$ commutes with the $\mu_r$-action.
Now using Proposition \ref{traceprop} and Theorem \ref{functisomthm},
\begin{equation*}
\begin{split}
\ch(U\circ V)|_{p_i(\zeta)\to\tr(\varphi^i\zeta,M)}
&=\tr(F_{U\circ V}(\varphi),F_{U\circ V}(M))\\
&=\tr(F_U(F_V^{(r)}(\varphi)),F_U(F_V^{(r)}(M)))\\
&=\ch(U)|_{p_i(\zeta)\to\tr(F_V^{(r)}(\varphi)^i\zeta,F_V^{(r)}(M))}\\
&=\ch(U)|_{p_i(\zeta)\to\tr(F_V(\varphi^i\zeta),F_V(M))}\\
&=\ch(U)|_{p_i(\zeta)\to\ch(V)|_{p_j\to\tr(\varphi^{ij}\zeta^j,M)}}\\
&=(\ch(U)\circ\ch(V))|_{p_i(\zeta)\to\tr(\varphi^i\zeta,M)},
\end{split}
\end{equation*}
as required.
\end{proof}
\begin{corollary} \label{specplethysmcor}
If $A$ is an $r$-species, and $B$ is a species with
$B(\emptyset)=\emptyset$, then $Z_{A\circ B}=Z_A\circ Z_B$.
\end{corollary}
\begin{proof}
Since $H^0(A\circ B)\cong H^0 A\circ H^0 B$, this is an immediate consequence.
\end{proof}

Note that all the results of this section would remain true if
$\mu_r$ were replaced by any finite group $G$ (and $\B_r$ by the category
of finite sets with a free $G$-action and bijections respecting
these actions, etc.) Thus there is a theory of ``$(G\wr\Ss)$-modules''
for all finite groups $G$. I do not know any uses for this
more general concept.
%%%%%%%%%%%%%%%%%%%%%%%%%%%%%%%%%%%%%%%%%%%%%%%%%%%%%%%%%%%%%%%%%%%%%%%%
\section{Stratifications and Trees}
We now begin the proof of Theorems \ref{ginzburgthm} and \ref{mainthm},
by analysing the nested-set stratifications defined in 
\cite{wonderful}. It is convenient to adopt a functorial
viewpoint, similar to that of $r$-species, but with the target category 
$\B_1$ replaced by the category $\Var_1$ of complex varieties and isomorphisms.
We will call a functor $\B_r\to\Var_1$ a \emph{$\B_r$-variety}.
Since a finite set can be regarded as a complex variety, this notion actually
includes that of $r$-species. If $A(r)$ is a $\B_r$-variety, we will use the
notation $A(r,n)$ for $A(r,[n]_r)$, which is a variety with an action of
$W(r,n)$. 

We first show that the varieties we have already introduced
fit into this pattern. For any object $I$ of $\B_r$, define
\[ V(r,I):=\left\{\begin{array}{cl}
\C^I/\{\text{constant fns}\},&\text{ if $r=1$, and}\\
\{(z_i)\in\C^I\,|\,z_{\zeta.i}=\zeta z_i,\ \forall i\in I,\zeta\in\mu_r\},
&\text{ if $r\geq 2$.}
\end{array}\right. \]
(If $I=\emptyset$, we interpret this to mean $\{0\}$ in either case.)
With the obvious definition on morphisms, this gives us a $\B_r$-variety
$V(r)$ (whose values happen to be
finite-dimensional vector spaces). We can identify $V(r,n):=V(r,[n]_r)$
with the vector space denoted $V(r,n)$ in the Introduction, 
by the map sending $(z_{(\zeta,j)})_{\zeta\in\mu_r,j\in[n]}$
to $(z_{(1,j)})_{j\in [n]}$.

Similarly, we define an $r$-species $\CalA(r)$ by letting
$\CalA(r,I)$ be the set of hyperplanes of $V(r,I)$ which are given
by an equation of the form $z_i=z_{i'}$ for some $i\neq i'$ in $I$.
It is easy to see that this coheres with the previous definition of
$\CalA(r,n)$; for instance, the hyperplanes
given by $z_i=z_{i'}$ where $i$ and $i'$ are in the same
$\mu_r$-orbit correspond to those in the previous definition
where one of the coordinates
is $0$. 

Finally, define $\B_r$-varieties $\M(r),\Mbar(r)$ by letting
$\M(r,I)$ be the complement in $\Pp(V(r,I))$ of the hyperplanes in 
$\CalA(r,I)$ (this is empty if $V(r,I)=\{0\}$),
and $\Mbar(r,I)$ its De Concini-Procesi compactification as defined
in the Introduction. For convenience in stating Theorem \ref{1stratthm} below,
we stipulate that if $|I|=1$,
$\Mbar(1,I)$ is actually a single point, not the empty set.

In this paper, a \emph{stratification} of a complex
variety $X$ will mean a quadruple $(A,\{X_\alpha\}_{\alpha\in A},
\{Y_\alpha\}_{\alpha\in A},\{\phi_\alpha\}_{\alpha\in A})$ such that:
\begin{enumerate}
\item $A$ is a finite set;
\item for each $\alpha\in A$, the \emph{stratum} $X_\alpha$ is a locally
closed subvariety of $X$;
\item $X$ is the disjoint union of all $X_\alpha$'s;
\item the closure of each $X_\alpha$ in $X$
is a union of strata;
\item for each $\alpha\in A$, $Y_\alpha$ is a complex variety
and $\phi_\alpha:X_\alpha\isomto Y_\alpha$ is an isomorphism.
\end{enumerate}
If $X$ is a $\B_r$-variety, a stratification of $X$ is defined in the same way,
but with all the data depending functorially on an object $I$ of $\B_r$.
Explicitly, it consists of an $r$-species $A$, a stratification
of $X(I)$ for all $I$ in which the indexing set is $A(I)$,
and collections of isomorphisms
\[ X(f)_{\alpha}:X_\alpha\isomto X_{A(f)(\alpha)},\
Y(f)_{\alpha}:Y_\alpha\isomto Y_{A(f)(\alpha)}, \]
for any morphism $f:I\isomto J$ in $\B_r$
and $\alpha\in A(I)$. These must satisfy the conditions that
$X(f):X(I)\isomto X(J)$ restricts to $X(f)_\alpha$ on $X_\alpha$, and
$\phi_{A(f)(\alpha)}\circ X(f)_\alpha=Y(f)_\alpha\circ\phi_\alpha$,
for all $f$ and $\alpha$. Given such a stratification of the $\B_r$-variety
$X$, we can define a $\B_r$-variety $Y$ by the rule
\[ Y(I)=\coprod_{\alpha\in A(I)} Y_\alpha, \]
using the isomorphisms $Y_{f,\alpha}$ to define $Y$ on morphisms.
We can regard $Y$ as being the $\B_r$-variety obtained
by ``cutting $X$ into pieces'' according to the stratifications.
We will say for short that we have a stratification
$X\strat Y$ of $\B_r$-varieties, leaving the other data understood.

Define a species $\T(1):\B_1\to\B_1$ by letting $\T(1,I)$
be the set of isomorphism classes of rooted trees, endowed with a bijection
between $I$ and the set of leaves, and satisfying the condition that
every internal (i.e.\ non-leaf) vertex has at least two edges in its fibre
(i.e.\ incident with the vertex and leading away from the root).
As the following pictures show, $|\T(1,0)|=0$, $|\T(1,1)|=1$, 
$|\T(1,2)|=1$, and $|\T(1,3)|=4$.

\medskip

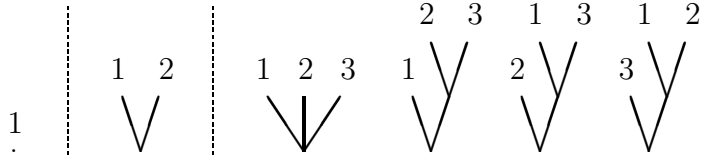
\setlength{\unitlength}{1500sp}%
\begingroup\makeatletter\ifx\SetFigFont\undefined%
\gdef\SetFigFont#1#2#3#4#5{%
  \reset@font\fontsize{#1}{#2pt}%
  \fontfamily{#3}\fontseries{#4}\fontshape{#5}%
  \selectfont}%
\fi\endgroup%
\begin{picture}(12766,2574)(768,-2223)
\thicklines
\put(4201,-2161){\line( 1, 3){300}}
\put(6901,-2161){\line(-2, 3){600}}
\put(6901,-2161){\line( 0, 1){900}}
\put(6901,-2161){\line( 2, 3){600}}
\put(9001,-2161){\line(-1, 3){300}}
\put(9001,-2161){\line( 1, 3){300}}
\put(9301,-1261){\line(-1, 3){300}}
\put(9301,-1261){\line( 1, 3){300}}
\put(10801,-2161){\line(-1, 3){300}}
\put(10801,-2161){\line( 1, 3){300}}
\put(11101,-1261){\line(-1, 3){300}}
\put(11101,-1261){\line( 1, 3){300}}
\put(12601,-2161){\line(-1, 3){300}}
\put(12601,-2161){\line( 1, 3){300}}
\put(12901,-1261){\line(-1, 3){300}}
\put(12901,-1261){\line( 1, 3){300}}
\put(4201,-2161){\line(-1, 3){300}}
\put(2101,-2161){\makebox(6.6667,10.0000){\SetFigFont{10}{12}{\rmdefault}{\mddefault}{\updefault}.}}
\put(13201,-61){\makebox(0,0)[lb]{\smash{\SetFigFont{12}{6.0}{\rmdefault}{\mddefault}{\updefault}2}}}
\thinlines
\multiput(3001,239)(0.00000,-120.00000){21}{\line( 0,-1){ 60.000}}
\multiput(5401,239)(0.00000,-120.00000){21}{\line( 0,-1){ 60.000}}
\put(2001,-1861){\makebox(0,0)[lb]{\smash{\SetFigFont{12}{6.0}{\rmdefault}{\mddefault}{\updefault}1}}}
\put(3701,-961){\makebox(0,0)[lb]{\smash{\SetFigFont{12}{6.0}{\rmdefault}{\mddefault}{\updefault}1}}}
\put(4501,-961){\makebox(0,0)[lb]{\smash{\SetFigFont{12}{6.0}{\rmdefault}{\mddefault}{\updefault}2}}}
\put(6101,-961){\makebox(0,0)[lb]{\smash{\SetFigFont{12}{6.0}{\rmdefault}{\mddefault}{\updefault}1}}}
\put(6801,-961){\makebox(0,0)[lb]{\smash{\SetFigFont{12}{6.0}{\rmdefault}{\mddefault}{\updefault}2}}}
\put(7501,-961){\makebox(0,0)[lb]{\smash{\SetFigFont{12}{6.0}{\rmdefault}{\mddefault}{\updefault}3}}}
\put(8501,-961){\makebox(0,0)[lb]{\smash{\SetFigFont{12}{6.0}{\rmdefault}{\mddefault}{\updefault}1}}}
\put(8801,-61){\makebox(0,0)[lb]{\smash{\SetFigFont{12}{6.0}{\rmdefault}{\mddefault}{\updefault}2}}}
\put(9601,-61){\makebox(0,0)[lb]{\smash{\SetFigFont{12}{6.0}{\rmdefault}{\mddefault}{\updefault}3}}}
\put(10301,-961){\makebox(0,0)[lb]{\smash{\SetFigFont{12}{6.0}{\rmdefault}{\mddefault}{\updefault}2}}}
\put(10601,-61){\makebox(0,0)[lb]{\smash{\SetFigFont{12}{6.0}{\rmdefault}{\mddefault}{\updefault}1}}}
\put(11401,-61){\makebox(0,0)[lb]{\smash{\SetFigFont{12}{6.0}{\rmdefault}{\mddefault}{\updefault}3}}}
\put(12101,-961){\makebox(0,0)[lb]{\smash{\SetFigFont{12}{6.0}{\rmdefault}{\mddefault}{\updefault}3}}}
\put(12401,-61){\makebox(0,0)[lb]{\smash{\SetFigFont{12}{6.0}{\rmdefault}{\mddefault}{\updefault}1}}}
\end{picture}

\medskip

\noindent
Abusing notation slightly,
we will usually regard $\T(1,I)$ as a set of trees representing the isomorphism
classes rather than the set of classes themselves. For $T\in\T(1,I)$, we will 
write $\Int(T)$ for the
set of internal vertices, and for every $v\in\Int(T)$, $\Fibre(v)$
for the fibre of $v$.

If $X$ is a $\B_1$-variety, define
a $\B_1$-variety $\mathbb{T}X$ by the rule
\[ \mathbb{T}X(I)=\coprod_{T\in\T(1,I)}\prod_{v\in\Int(T)} X(\Fibre(v)). \]
When $|I|=1$, this is an empty product, i.e.\ a single point. To illustrate
this definition, suppose $I=[3]$.
The terms of the disjoint union corresponding to the four elements of
$\T(1,3)$ are isomorphic to $X[3]$, $X[2]\times X[2]$, $X[2]\times X[2]$,
and $X[2]\times X[2]$ respectively. But the key point is how these
varieties are treated functorially: for instance, the transposition 
$(1\ 2)\in S_3$ permutes the second and third varieties, because it 
interchanges the corresponding trees; it fixes the fourth variety,
acting on the $X[2]$ factor coresponding to the root as the identity,
and on the other $X[2]$ factor as $(1\ 2)$, since that is how the
corresponding automorphism of this tree acts on the fibres of the root and
of the other vertex.

The following result is well known in the moduli space context
(see \cite[Chapter III, \S2.8]{maninbook}).
It is also deducible as a special case of the results of 
\cite{wonderful}, in a similar way to the next Theorem.
\begin{theorem} \label{1stratthm}
We have a stratification $\Mbar(1)\strat\mathbb{T}\M(1)$ of $\B_1$-varieties.
\end{theorem}
\noindent
This Theorem says in part that each variety $\Mbar(1,n)$ has a stratification
where the strata are indexed by $\T(1,n)$, and the stratum corresponding to
a tree $T$ is isomorphic to $\prod_{v\in\Int(T)}\M(1,\Fibre(v))$.
Crucially, it also says that these strata are permuted and acted upon by
$S_n$ according to its action on these trees. 

Intuitively, one can think of $\Mbar(1,n)$ as the set of ways in which
$n$ distinct points in the complex plane 
can ``collapse to the origin'': the vertices
of a tree in $\T(1,n)$, read from the leaves to the root, record the
order in which the points coalesce, and the varieties $\M(1,\Fibre(v))$
record the relative configuration of each set of points ``just before''
they coincide.

Now suppose $r\geq 2$. Define an $r$-species $\T(r):\B_r\to\B_1$
by letting $\T(r,I)$ be the subset of $\T(1,I)$ consisting of
(isomorphism classes of) trees $T$ satisfying the following extra
conditions:
\begin{enumerate}
\item $T$ is $\mu_r$-stable, i.e.\ the action of $\mu_r$ on the leaves
of $T$ extends to an action of $\mu_r$ on $T$ by rooted tree automorphisms;
\item every $\mu_r$-fixed vertex of $T$ has at most one $\mu_r$-fixed
edge in its fibre;
\item for every $\mu_r$-fixed vertex, $\mu_r$ acts freely on the set
of edges in the fibre which are not $\mu_r$-fixed.
\end{enumerate} 
Note that condition (2) implies that the
$\mu_r$-fixed vertices lie on a single unbranched path from the root;
condition (3) then implies that $\mu_r$ acts freely on the
set of non-$\mu_r$-fixed vertices. 
For example, the following pictures show that
$|\T(2,1)|=1$ and $|\T(2,2)|=5$ (we have written the 
elements of $[n]_2$ as $\pm i$ rather than $(\pm 1,i)$).

\medskip

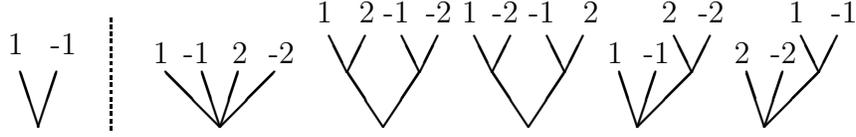
\setlength{\unitlength}{1500sp}%
\begingroup\makeatletter\ifx\SetFigFont\undefined%
\gdef\SetFigFont#1#2#3#4#5{%
  \reset@font\fontsize{#1}{#2pt}%
  \fontfamily{#3}\fontseries{#4}\fontshape{#5}%
  \selectfont}%
\fi\endgroup%
\begin{picture}(14600,2172)(0,-1373)
\thicklines
\put(901,-1261){\line(-1, 3){300}}
\put(901,-1261){\line( 1, 3){300}}
\put(3901,-1261){\line(-1, 1){900}}
\put(3901,-1261){\line(-1, 3){300}}
\put(3901,-1261){\line( 1, 3){300}}
\put(3901,-1261){\line( 1, 1){900}}
\put(6601,-1261){\line(-2, 3){600}}
\put(6601,-1261){\line( 2, 3){600}}
\put(6001,-361){\line(-1, 2){300}}
\put(6001,-361){\line( 1, 2){300}}
\put(7201,-361){\line(-1, 2){300}}
\put(7201,-361){\line( 1, 2){300}}
\put(9001,-1261){\line(-2, 3){600}}
\put(8401,-361){\line(-1, 2){300}}
\put(8401,-361){\line( 1, 2){300}}
\put(9001,-1261){\line( 2, 3){600}}
\put(9601,-361){\line(-1, 2){300}}
\put(9601,-361){\line( 1, 2){300}}
\put(10801,-1261){\line(-1, 3){300}}
\put(10801,-1261){\line( 1, 3){300}}
\put(10801,-1261){\line( 1, 1){900}}
\put(11701,-361){\line(-1, 2){300}}
\put(11701,-361){\line( 1, 2){300}}
\put(12901,-1261){\line(-1, 3){300}}
\put(12901,-1261){\line( 1, 3){300}}
\thinlines
%\put(901,-1261){\line(-1, 3){300}}
\multiput(2101,539)(0.00000,-120.00000){16}{\line( 0,-1){ 60.000}}
\put(14001,464){\makebox(0,0)[lb]{\smash{\SetFigFont{12}{14.4}{\rmdefault}{\mddefault}{\updefault}-1}}}
\thicklines
\put(12901,-1261){\line( 1, 1){900}}
\put(13801,-361){\line(-1, 2){300}}
\put(13801,-361){\line( 1, 2){300}}
\put(401,-61){\makebox(0,0)[lb]{\smash{\SetFigFont{12}{14.4}{\rmdefault}{\mddefault}{\updefault}1}}}
\put(1101,-61){\makebox(0,0)[lb]{\smash{\SetFigFont{12}{14.4}{\rmdefault}{\mddefault}{\updefault}-1}}}
\put(2801,-211){\makebox(0,0)[lb]{\smash{\SetFigFont{12}{14.4}{\rmdefault}{\mddefault}{\updefault}1}}}
\put(3301,-211){\makebox(0,0)[lb]{\smash{\SetFigFont{12}{14.4}{\rmdefault}{\mddefault}{\updefault}-1}}}
\put(4101,-211){\makebox(0,0)[lb]{\smash{\SetFigFont{12}{14.4}{\rmdefault}{\mddefault}{\updefault}2}}}
\put(4701,-211){\makebox(0,0)[lb]{\smash{\SetFigFont{12}{14.4}{\rmdefault}{\mddefault}{\updefault}-2}}}
\put(5501,464){\makebox(0,0)[lb]{\smash{\SetFigFont{12}{14.4}{\rmdefault}{\mddefault}{\updefault}1}}}
\put(6201,464){\makebox(0,0)[lb]{\smash{\SetFigFont{12}{14.4}{\rmdefault}{\mddefault}{\updefault}2}}}
\put(6601,464){\makebox(0,0)[lb]{\smash{\SetFigFont{12}{14.4}{\rmdefault}{\mddefault}{\updefault}-1}}}
\put(7301,464){\makebox(0,0)[lb]{\smash{\SetFigFont{12}{14.4}{\rmdefault}{\mddefault}{\updefault}-2}}}
\put(7901,464){\makebox(0,0)[lb]{\smash{\SetFigFont{12}{14.4}{\rmdefault}{\mddefault}{\updefault}1}}}
\put(8401,464){\makebox(0,0)[lb]{\smash{\SetFigFont{12}{14.4}{\rmdefault}{\mddefault}{\updefault}-2}}}
\put(9001,464){\makebox(0,0)[lb]{\smash{\SetFigFont{12}{14.4}{\rmdefault}{\mddefault}{\updefault}-1}}}
\put(9901,464){\makebox(0,0)[lb]{\smash{\SetFigFont{12}{14.4}{\rmdefault}{\mddefault}{\updefault}2}}}
\put(10301,-211){\makebox(0,0)[lb]{\smash{\SetFigFont{12}{14.4}{\rmdefault}{\mddefault}{\updefault}1}}}
\put(10901,-211){\makebox(0,0)[lb]{\smash{\SetFigFont{12}{14.4}{\rmdefault}{\mddefault}{\updefault}-1}}}
\put(11201,464){\makebox(0,0)[lb]{\smash{\SetFigFont{12}{14.4}{\rmdefault}{\mddefault}{\updefault}2}}}
\put(11801,464){\makebox(0,0)[lb]{\smash{\SetFigFont{12}{14.4}{\rmdefault}{\mddefault}{\updefault}-2}}}
\put(12401,-211){\makebox(0,0)[lb]{\smash{\SetFigFont{12}{14.4}{\rmdefault}{\mddefault}{\updefault}2}}}
\put(13001,-211){\makebox(0,0)[lb]{\smash{\SetFigFont{12}{14.4}{\rmdefault}{\mddefault}{\updefault}-2}}}
\put(13301,464){\makebox(0,0)[lb]{\smash{\SetFigFont{12}{14.4}{\rmdefault}{\mddefault}{\updefault}1}}}
\end{picture}

\medskip

\noindent
Here are three of the $47$ elements of $\T(2,3)$:

\medskip

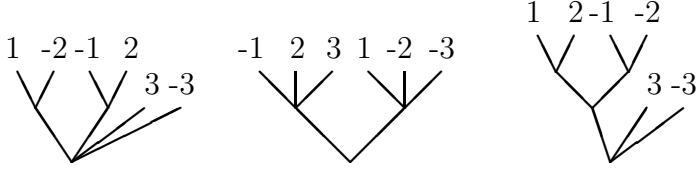
\setlength{\unitlength}{1500sp}%
\begingroup\makeatletter\ifx\SetFigFont\undefined%
\gdef\SetFigFont#1#2#3#4#5{%
  \reset@font\fontsize{#1}{#2pt}%
  \fontfamily{#3}\fontseries{#4}\fontshape{#5}%
  \selectfont}%
\fi\endgroup%
\begin{picture}(12000,2493)(-500,-2194)
\thicklines
\put(1501,-2161){\line( 2, 3){600}}
\put(1501,-2161){\line( 4, 3){1200}}
\put(901,-1261){\line(-1, 2){300}}
\put(901,-1261){\line( 1, 2){300}}
\put(2101,-1261){\line(-1, 2){300}}
\put(2101,-1261){\line( 1, 2){300}}
\put(1501,-2161){\line( 2, 1){1800}}
\put(6101,-2161){\line(-1, 1){900}}
\put(6101,-2161){\line( 1, 1){900}}
\put(5201,-1261){\line(-1, 1){600}}
\put(5201,-1261){\line( 0, 1){600}}
\put(5201,-1261){\line( 1, 1){600}}
\put(7001,-1261){\line(-1, 1){600}}
\put(7001,-1261){\line( 0, 1){600}}
\put(7001,-1261){\line( 1, 1){600}}
\put(10401,-2161){\line(-1, 3){300}}
\put(10101,-1261){\line(-1, 1){600}}
\put(10101,-1261){\line( 1, 1){600}}
\put(9501,-661){\line(-1, 2){300}}
\put(9501,-661){\line( 1, 2){300}}
\put(1501,-2161){\line(-2, 3){600}}
\put(10701,-661){\line(-1, 2){300}}
\put(11401,-1036){\makebox(0,0)[lb]{\smash{\SetFigFont{12}{14.4}{\rmdefault}{\mddefault}{\updefault}-3}}}
\put(10701,-661){\line( 1, 2){300}}
\put(10401,-2161){\line( 2, 3){600}}
\put(10401,-2161){\line( 4, 3){1200}}
\put(401,-436){\makebox(0,0)[lb]{\smash{\SetFigFont{12}{14.4}{\rmdefault}{\mddefault}{\updefault}1}}}
\put(1001,-436){\makebox(0,0)[lb]{\smash{\SetFigFont{12}{14.4}{\rmdefault}{\mddefault}{\updefault}-2}}}
\put(1551,-436){\makebox(0,0)[lb]{\smash{\SetFigFont{12}{14.4}{\rmdefault}{\mddefault}{\updefault}-1}}}
\put(2351,-436){\makebox(0,0)[lb]{\smash{\SetFigFont{12}{14.4}{\rmdefault}{\mddefault}{\updefault}2}}}
\put(2701,-1036){\makebox(0,0)[lb]{\smash{\SetFigFont{12}{14.4}{\rmdefault}{\mddefault}{\updefault}3}}}
\put(3101,-1036){\makebox(0,0)[lb]{\smash{\SetFigFont{12}{14.4}{\rmdefault}{\mddefault}{\updefault}-3}}}
\put(4251,-436){\makebox(0,0)[lb]{\smash{\SetFigFont{12}{14.4}{\rmdefault}{\mddefault}{\updefault}-1}}}
\put(5101,-436){\makebox(0,0)[lb]{\smash{\SetFigFont{12}{14.4}{\rmdefault}{\mddefault}{\updefault}2}}}
\put(5701,-436){\makebox(0,0)[lb]{\smash{\SetFigFont{12}{14.4}{\rmdefault}{\mddefault}{\updefault}3}}}
\put(6201,-436){\makebox(0,0)[lb]{\smash{\SetFigFont{12}{14.4}{\rmdefault}{\mddefault}{\updefault}1}}}
\put(6701,-436){\makebox(0,0)[lb]{\smash{\SetFigFont{12}{14.4}{\rmdefault}{\mddefault}{\updefault}-2}}}
\put(7401,-436){\makebox(0,0)[lb]{\smash{\SetFigFont{12}{14.4}{\rmdefault}{\mddefault}{\updefault}-3}}}
\put(9001,164){\makebox(0,0)[lb]{\smash{\SetFigFont{12}{14.4}{\rmdefault}{\mddefault}{\updefault}1}}}
\put(9701,164){\makebox(0,0)[lb]{\smash{\SetFigFont{12}{14.4}{\rmdefault}{\mddefault}{\updefault}2}}}
\put(10051,164){\makebox(0,0)[lb]{\smash{\SetFigFont{12}{14.4}{\rmdefault}{\mddefault}{\updefault}-1}}}
\put(10801,164){\makebox(0,0)[lb]{\smash{\SetFigFont{12}{14.4}{\rmdefault}{\mddefault}{\updefault}-2}}}
\put(11001,-1036){\makebox(0,0)[lb]{\smash{\SetFigFont{12}{14.4}{\rmdefault}{\mddefault}{\updefault}3}}}
\end{picture}

\medskip

\noindent
For $T\in\T(r,I)$, write $\Int(T)^{\mu_r}$ for the set of
$\mu_r$-fixed vertices of $T$ (which must all be internal). If 
$v\in\Int(T)^{\mu_r}$, let $\Fibre(v)^\circ$ be the set of edges in
the fibre which are not $\mu_r$-fixed; by (3), this is an object
of $\B_r$. Write $\Orb(T)$ for the set of $\mu_r$-orbits in the 
set of non-$\mu_r$-fixed
vertices in $\Int(T)$. 

If $X$ is a $\B_r$-variety and
$Y$ is a $\B_1$-variety,
define a $\B_r$-variety $\mathbb{T}_r(X,Y)$ by the rule
\[ \mathbb{T}_r(X,Y)(I)=\negthickspace\coprod_{T\in\T(r,I)}
\left(
\prod_{v\in\Int(T)^{\mu_r}}
\negthickspace
X(\Fibre(v)^\circ)\ \times\negthickspace\negthickspace
\prod_{\CalO\in\Orb(T)}\negthickspace
Y(\Fibre(\CalO))\right), \]
where $Y(\Fibre(\CalO))$ denotes the
limit of all $Y(\Fibre(v))$ for $v\in\CalO$: in other words,
\[ \{(y_v)\in\prod_{v\in\CalO} Y(\Fibre(v))\,|\,
y_{\zeta.v}=Y(\alpha_{\zeta,v})(y_v),
\ \forall \zeta\in\mu_r,v\in\CalO\}, \]
where $\alpha_{\zeta,v}:\Fibre(v)\isomto\Fibre(\zeta.v)$ is the
bijection induced by the $\mu_r$-action on $T$.
\begin{theorem} \label{rstratthm}
We have a stratification 
$\Mbar(r)\strat\mathbb{T}_r(\M(r),\M(1))$ of $\B_r$-varieties.
\end{theorem}
\noindent
As in the $r=1$ case, one can think intuitively of $\Mbar(r,n)$
as the set of ways of collapsing $n$ points in $\C$ to the origin,
except that now we care about the $\mu_r$-orbits of the points,
not just the points themselves. For example, the following 
is the sequence of events
in the stratum corresponding to the first element of $\T(2,3)$ depicted
above.
\begin{enumerate}
\item The first point and the negative of the second point coalesce,
which forces their negatives to coalesce also; the configuration of either
coalescing pair is recorded as a point (rather, the unique point)
of $\M(1,2)$.
\item The two new combined points and the original third point and its negative
coalesce, i.e.\ all become zero; the configuration of these two $\mu_2$-orbits
is recorded as a point of $\M(2,2)$.
\end{enumerate}
In general, the $\mu_r$-fixed vertices represent the coalescing of points at
the origin, and the non-$\mu_r$-fixed vertices represent
the coalescing of points away from the origin (occurring simultaneously
in $\mu_r$-orbits); this explains conditions (2) and (3) in the definition
of $\T(r,I)$.

Theorem \ref{rstratthm} is deduced from \cite{wonderful} as follows.
\begin{proof}
First consider the general situation of an essential irreducible
arrangement $\CalA$ in a vector space $V$. Let $\M$, $\Mbar$, $\F$ be
as in the Introduction.
The results of \cite[\S\!\S 4.2, 4.3]{wonderful} are stated in terms
of the building set $\mathcal{G}$, which in our case is the set
of subspaces of $V^*$ perpendicular to the elements of $\F$. It
is easy to translate them into the following statements about $\F$:
\begin{enumerate}
\item  The irreducible components
of $\Mbar\setminus\M$ are in bijection with the set
$\F^*:=\F\setminus\{\{0\}\}$. If $\overline{D}_X$ denotes the component
corresponding to $X\in\F^*$, then the set of points in $\Mbar(r,I)$
whose image in $\Pp(V)$ lies in $\Pp(X)$ is precisely
$\bigcup_{Y\in\F^*, Y\subseteq X}\overline{D}_Y$.
\item  We call a subset $\mathcal{S}\subseteq\F$ \emph{nested}
if for all collections $\{X_1,\cdots,X_s\}$, $s\geq 2$, of pairwise
incomparable elements of $\mathcal{S}$, $\bigcap_{i=1}^{s} X_i\not\in\F$.
(The empty subset and singleton sets are automatically nested.) If
$\mathcal{S}\subseteq\F^*$, then
\[ \bigcap_{X\in\mathcal{S}} \overline{D}_X\neq\emptyset
\Longleftrightarrow \text{ $\mathcal{S}$ is nested.} \]
\item  Hence we have a disjoint union
\[ \Mbar=\bigcup_{\substack{\mathcal{S}\subseteq\F^*\\\text{$\mathcal{S}$
nested}}} \M_{\mathcal{S}}, \]
where $\M_{\mathcal{S}}=\{x\in\Mbar\,|\,
\mathcal{S}=\{X\in\F^*\,|\,x\in\overline{D}_X\}\}$. Note that $\M_\emptyset
=\M$.
\item  If $\mathcal{S}\subseteq\F^*$ is nested, 
$\M_{\mathcal{S}}\cong\prod_{X\in\mathcal{S}\cup\{\{0\}\}}\M_{\mathcal{S},X}$, 
where $\M_{\mathcal{S},X}$ is the projectivized
hyperplane complement of the arrangement $\CalA_{\mathcal{S},X}$
in $(\bigcap_{Y\in\mathcal{S},
Y\supset X} Y)/X$ induced by those hyperplanes in $\CalA$ which
contain $X$.
\end{enumerate}
(This last point requires a little thought, since 
\cite[Theorem 4.3]{wonderful} is stated in terms of
$\overline{\M_{\mathcal{S}}}=\bigcap_{X\in\mathcal{S}}\overline{D}_X$
instead of $\M_{\mathcal{S}}$.)

To prove the Theorem, we must unravel
what all this means in the case $\CalA=\CalA(r,I)$. Most of the details
will be left to the reader.
If $J$ is an $\mu_r$-stable subset of $I$,
let 
\[ X_{J}:=\{(z_i)\in V(r,I)\,|\, z_j=0,\ \forall j\in J\}. \]
If in addition $\pi=\{J_1,\cdots,J_r\}$ is a partition of $J$ on which
$\mu_r$ acts freely and transitively, let
\[ X_{J,\pi}=\{(z_i)\in V(r,I)\,|\, z_j=z_{j'},\ \forall j,j'\in J_k, 
1\leq k\leq r\}. \]
It is not hard to see that
\[ \F=\{X_{J}\,|\,J\neq\emptyset\}
\cup\{X_{J,\pi}\,|\,|J|\geq 2r\}. \]
(Note that $\{0\}=X_I$.) It is also straightforward
to show that
$\mathcal{S}\subseteq\F^*$ is nested if and only if the following 
conditions hold:
\begin{enumerate}
\item  the subsets $J\subset I$ such that $X_J\in\mathcal{S}$
are totally ordered;
\item  if $X_J,X_{K,\pi}\in\mathcal{S}$, 
then either $J\cap K=\emptyset$ or $J\supseteq K$;
\item  if $X_{J,\pi},
X_{K,\pi'}\in\mathcal{S}$, then either $J\cap K=\emptyset$ or
one of $X_{J,\pi}$, $X_{K,\pi'}$ contains the other.
\end{enumerate}
(The $r=2$ case is contained in \cite[\S4]{yuzvinsky}.)

Now let $T\in\T(r,I)$. For any $v\in\Int(T)$, let $J_v\subseteq I$
be the subset corresponding to the leaves of the subtree with root $v$.
If $v\in\Int(T)^{\mu_r}$, then $J_v$ is $\mu_r$-stable (if $v$ is the root,
$J_v=I$). If $\CalO\in\Orb(T)$, then $J_\CalO:=\bigcup_{v\in\CalO}J_v$
is $\mu_r$-stable and $\pi_\CalO:=\{J_v\,|\,v\in\CalO\}$ is a partition
of $J_\CalO$ on which $\mu_r$ acts freely and transitively. It is easy to
show that
\[ \mathcal{S}_\T:=\{X_{J_v}\,|\,v\in\Int(T)^{\mu_r}\text{ not the root}\}
\cup\{X_{J_\CalO,\pi_\CalO}\,|\,\CalO\in\Orb(T)\} \]
is nested, and that this gives a bijection
between $\T(r,I)$ and the set of nested subsets of $\F^*$. (The tree
in which the root is the sole internal vertex corresponds to the
empty set.) Writing $\M_T$ for $\M_{\mathcal{S}_\T}$ gives us a
disjoint union
\begin{equation} \label{firststrateqn}
\Mbar(r,I)=\bigcup_{T\in\T(r,I)}\M_T. 
\end{equation}
Now we express $\M_T$ as a product over $\Int(T)^{\mu_r}
\amalg\Orb(T)$ using property (4) above. If $v\in\Int(T)^{\mu_r}$,
then the elements $Y\in\mathcal{S}_T$ such that $Y\supseteq X_{J_v}$
are precisely the $X_{J_{v'}}$ and $X_{J_{\CalO},\pi_\CalO}$ such that
the vertex $v'$ or the orbit $\CalO$ lies in the subtree with root $v$.
Among these, the minimal subspaces containing $X_{J_v}$ strictly
are those corresponding to $v'$
or $\CalO$ lying in the fibre of $v$. Hence one gets that the arrangement
$\CalA_{\mathcal{S}_T,X_{J_v}}$ is of type $\CalA(r,\Fibre(v)^\circ)$.
Similarly, for $\CalO\in\Orb(T)$ the arrangement 
$\CalA_{\mathcal{S}_T,X_{J_\CalO,\pi_\CalO}}$ is of type 
$\CalA(1,\Fibre(\CalO))$. Hence
\begin{equation} \label{secondstrateqn}
\M_T\cong\prod_{v\in\Int(T)^{\mu_r}}
\negthickspace
\M(r,\Fibre(v)^\circ)\ \times
\prod_{\CalO\in\Orb(T)}\negthickspace
\M(1,\Fibre(\CalO)).
\end{equation}
Finally, one can verify that \eqref{firststrateqn} and \eqref{secondstrateqn}
constitute a stratification of $\B_r$-varieties.
\end{proof}

Theorems \ref{1stratthm} and \ref{rstratthm} motivate us to examine
the operations $\mathbb{T}$ and $\mathbb{T}_r$ more closely. In
particular, we want to relate them to the operations of addition,
multiplication, and substitution for $\B_r$-varieties, which are defined
in exactly the same way as for $r$-species (using the usual disjoint union
and product of varieties). We will thus obtain equations which capture 
the recursive nature of these trees.

The following result is one of many variants of
\cite[\S3.1, Theorem 2]{bll}, a key step in combinatorial Lagrange Inversion.
\begin{proposition} \label{1treeprop}
We have the following isomorphism of species:
\[ \T(1)\cong E(1)_1 + (E(1)_{\geq 2}\circ \T(1)). \]
\end{proposition}
\begin{proof}
The first term on the right-hand side
takes care of the case that $|I|=1$. Otherwise, any tree $T\in\T(1,I)$ 
can be split at the root to create an assembly of two or more smaller trees:
more formally, there is a unique partition $\pi_T\in\Par(I)$ and unique
$T_J\in\T(1,J)$
for all $J\in\pi_T$, such that $T$ is obtained by joining all the roots of
the trees $T_J$ to a new vertex which then becomes the root.
By definition of substitution of species, this is exactly 
encapsulated in the second term on the right-hand side.
\end{proof}
\noindent
Applying Corollary \ref{specplethysmcor} (or rather the $r=1$ case of it,
which is Joyal's result), we deduce:
\begin{equation} \label{1pletheqn}
Z_{\T(1)}=p_1+Z_{E(1)_{\geq 2}}\circ Z_{\T(1)}.
\end{equation}
In other words, $Z_{\T(1)}$ is the plethystic inverse of 
$p_1-Z_{E(1)_{\geq 2}}$.
As with Theorem \ref{ginzburgthm}, \eqref{1pletheqn} gives
a recursive formula for the terms of $Z_{\T(1)}$.

Proposition \ref{1treeprop} has an ``enriched'' version, implicit 
in \cite{ginzburgkapranov}.
\begin{proposition} \label{1plethprop}
If $X$ is a $\B_1$-variety with $X(I)=\emptyset$ for $|I|\leq 1$,
we have an isomorphism of $\B_1$-varieties
\[ \mathbb{T}X\cong E(1)_1 + (X\circ\mathbb{T}X). \]
\end{proposition}
\begin{proof}
The proof is essentially the same as that of Proposition \ref{1treeprop}, 
but this time we have to keep track of the varieties $X(\Fibre(v))$
attached to each $v\in\Int(T)$. Again, the first term on the right-hand side
takes care of the case that $|I|=1$. Otherwise, if
$\pi_T$ and $T_J$ are defined as before, it is clear that
\[ \prod_{v\in\Int(T)}
X(\Fibre(v))\cong X(\pi_T)\times
\prod_{J\in\pi_T}\prod_{v\in\Int(T_J)}
X(\Fibre(v)). \]
So the result follows by definition of substitution.
\end{proof}
\noindent
All that remains in order to prove Theorem \ref{ginzburgthm} is to translate
this isomorphism of $\B_r$-varieties into an isomorphism of their
cohomologies in a suitable category, and apply the appropriate variant
of Theorem \ref{plethysmthm}. This will be done in the next two sections.
For now, note that if $X$ in Proposition \ref{1plethprop} happens to be
a species (i.e.\ all the varieties are finite sets), then so is
$\mathbb{T}X$, and we deduce from Proposition \ref{1plethprop}
that $Z_{\mathbb{T}X}$ is the plethystic inverse of $p_1-Z_X$.

The $r\geq 2$ case is slightly different.
\begin{proposition} \label{rtreeprop}
For $r\geq 2$, we have the following isomorphism of $r$-species:
\[ \T(r)\cong (E(r)_0 + \T(r))\cdot
(E(r)_+\circ \T(1)). \]
\end{proposition}
\begin{proof}
Define a sub-$r$-species $\U(r)$ of $\T(r)$, by letting $\U(r,I)$
be the set of (isomorphism classes of) trees in $\T(r,I)$ in which
the root is the sole $\mu_r$-fixed vertex. Clearly any $T\in\T(r,I)$
which is not in $\U(r,I)$
determines a decomposition $(I_1,I_2)\in\Decomp(I)^{\mu_r}$, a
tree $T_1\in\T(r,I_1)$, and a tree $T_2\in\U(r,I_2)$, such that $T$
is obtained by joining the root of $T_1$ to the root of $T_2$ and making
the latter the root of the resulting tree (this new edge then becomes
the unique $\mu_r$-fixed edge in the fibre of the root). For $T\in\U(r,I)$,
we can think of the decomposition having empty $I_1$. Hence
\begin{equation*}
\T(r)\cong (E(r)_0 + \T(r))\cdot \U(r).
\end{equation*} 
But if $T\in\U(r,I)$, the partition $\pi_T$ defined in the proof
of Proposition \ref{1treeprop} clearly lies in
$\Par(I)_{\B_r}$, and the subtrees
$T_J\in\T(1,J)$ for each $J\in\pi_T$ satisfy
$T_{J}\cong T_{\zeta.J}$ via an isomorphism which
acts on the leaves by the action of $\zeta$. So by definition of substitution,
\begin{equation*}
\U(r)\cong E(r)_+\circ \T(1),
\end{equation*}
whence the result.
\end{proof}
\noindent
Applying Corollary \ref{specplethysmcor}, we deduce:
\begin{equation}
Z_{\T(r)}=(1+Z_{\T(r)})(Z_{E(r)_+}\circ Z_{\T(1)}),
\end{equation}
which means that $1+Z_{\T(r)}$ is the inverse
of $1-Z_{E(r)_+}\circ Z_{\T(1)}$ in $\Aa(r)$.

The ``enriched'' version is as follows:
\begin{proposition} \label{rplethprop}
If $X$ is a $\B_r$-variety such that $X(\emptyset)=\emptyset$, and $Y$
is a $\B_1$-variety such that $Y(I)=\emptyset$ if $|I|\leq 1$,
we have an isomorphism of $\B_r$-varieties
\[ \mathbb{T}_r (X,Y)\cong (E(r)_0 + \mathbb{T}_r (X,Y))\cdot
(X\circ \mathbb{T}Y). \] 
\end{proposition}
\begin{proof}
As with Proposition \ref{1plethprop}, the proof here uses the same
combinatorial analysis as in Proposition \ref{rtreeprop}. 
Two points deserve comment:
firstly, if $T\in\T(r,I)\setminus\U(r,I)$, the factors in the
$T$-term of $\mathbb{T}_r (X,Y)(I)$
corresponding to vertices in the subtree $T_1$ are precisely the same
as the $T_1$-term of $\mathbb{T}_r (X,Y)(I_1)$; and secondly, if $T\in\U(r,I)$,
the $T$-term of $\mathbb{T}_r (A,B)(I)$ is exactly
$X(\pi_T)\times\prod_{\CalO\in\mu_r\setminus\pi_T} \mathbb{T}Y(\CalO)$.
\end{proof}
\noindent
Again, if $X$ is in fact an $r$-species and $Y$ is a species,
we deduce that $1+Z_{\mathbb{T}_r(X,Y)}$ is the inverse of
$1-Z_X\circ Z_{\mathbb{T}Y}$. As for the case we are actually interested in,
where $X=\M(r)$ and $Y=\M(1)$, that requires a more technically sophisticated
approach.
%%%%%%%%%%%%%%%%%%%%%%%%%%%%%%%%%%%%%%%%%%%%%%%%%%%%%%%%%%%%%%%%%%%%%
\section{Bigraded $\B_r$-modules}
In order to apply the theory of $\B_r$-modules 
to our representations on cohomology of hyperplane complements and their
compactifications,
we will need a variant of Theorem \ref{plethysmthm} for bigraded
$\B_r$-modules. A \emph{bigrading} on a vector space $M$ means an
$(\N\times\N)$-grading, i.e.\ a direct sum decomposition $M=\oplus_{a,b\geq 0}
M_{a,b}$. Let $\C$-$\mathbf{mod}^\mathrm{bigr}$ be the category of bigraded
finite-dimensional vector spaces, and
$\C\mu_r$-$\mathbf{mod}^\mathrm{bigr}$ the category of bigraded representations
of $\mu_r$. A 
\emph{bigraded $\B_r$-module} is a functor 
$U:\B_r\to\C$-$\mathbf{mod}^\mathrm{bigr}$; thus we have
direct sum decompositions
$U(I)=\oplus_{a,b\geq 0} U(I)_{a,b}$ for all $I\in\B_r$, respected by
the maps $U(f)$ for $f:I\isomto J$. Any ungraded vector space or $\B_r$-module
can be regarded as a bigraded object, concentrated in bidegree $(0,0)$.

If $M\in\C$-$\mathbf{mod}^\mathrm{bigr}$ and $\psi\in\End(M)$ respects the
bigrading, define the \emph{bigraded trace}
\[ \tr_{q,t}(\psi,M)=\sum_{a,b\geq 0}\tr(\psi|_{M_{a,b}},M_{a,b})\, 
q^a (-t)^b 
\in\C[q,t]. \]
If $U$ is a bigraded $\B_r$-module,
define its \emph{bigraded characteristic} $\ch_{q,t}(U)\in\Aa(r)[q,t]
:=\Aa(r)\otimes\C[q,t]$ by
\begin{equation*}
\begin{split} 
\ch_{q,t}(U) &:=\sum_{n\geq 0} \sum_{a,b\geq 0} \ch_{W(r,n)}(U([n]_r)_{a,b})
\, q^a (-t)^b\\
&=\sum_{n\geq 0}\frac{1}{r^n n!}\sum_{w\in W(r,n)}\tr_{q,t}(w,U([n]_r))\,p_w.
\end{split}
\end{equation*}
We also write $\ch_q(U)$ for the specialization $\ch_{q,t}(U)|_{t\to 1}$.
Note that sending $t\to 1$ does not entirely forget the second grading,
since there is still the sign $(-1)^b$.

Any bounded bigraded $\B_r$-module $U$ gives rise
to a polynomial functor $F_U^{\mathrm{bigr}}:
\C\mu_r$-$\mathbf{mod}^\mathrm{bigr}\to\C$-$\mathbf{mod}^\mathrm{bigr}$,
defined on objects by the same formula as before:
\[  F_U^{\mathrm{bigr}}(M):=\bigoplus_{n\geq 0}
(U([n]_r)\otimes M^{\otimes n})^{W(r,n)}. \]
The bigrading of the tensor product of bigraded vector spaces
is defined in the usual way. The novelty we want to introduce is in the
action of $W(r,n)$ on $M^{\otimes n}$; this should again be given by the
permutation of the tensor factors composed with the actions of $\mu_r$ on
each, but the permutation should be accompanied by the graded-commutative sign
convention with respect to the second of the two gradings. 
In short, the $n$th tensor power should be taken in the symmetric monoidal 
category $g\mathrm{Vect}^-$ of \cite[(1.3.18)]{ginzburgkapranov};
this is the usual category of graded vector spaces, but 
with a commutativity isomorphism which differs by a sign:
\[ v\otimes w\mapsto (-1)^{\deg(v)\deg(w)} w\otimes v, \]
where in our case $\deg$ refers to the second grading.
What this means more explicitly is that
if $m_s\in M_{a_s,b_s}$ for all $1\leq s\leq n$, and $\tilde{w}\in S_n$, then
\[ \tilde{w}.(m_1\otimes\cdots\otimes m_n)
:=\varepsilon(\tilde{w},(b_i))\,
m_{\tilde{w}^{-1}(1)}\otimes\cdots\otimes m_{\tilde{w}^{-1}(n)}, \]
where $\varepsilon:S_n\times\N^n\to\{\pm 1\}$ is uniquely defined by
\begin{enumerate}
\item $\varepsilon(\tilde{y}\tilde{w},(b_i))=
\varepsilon(\tilde{y},\tilde{w}.(b_i))\varepsilon(\tilde{w},(b_i))$
for all $\tilde{y},\tilde{w}\in S_n$, $(b_i)\in\N^n$, and
\item $\varepsilon((k\ k+1),(b_i))=(-1)^{b_k b_{k+1}}$.
\end{enumerate}
The analogue of Proposition \ref{traceprop} is:
\begin{proposition} \label{gtraceprop}
If $U$ is a bounded
bigraded $\B_r$-module, $M\in\C\mu_r$-$\mathbf{mod}^\mathrm{bigr}$, and 
$\varphi\in\End(M)$ commutes with the $\mu_r$-action and
respects the bigrading, then
\[ \ch_{q,t}(U)|_{p_i(\zeta)\to\tr_{q,t}(\varphi^i\zeta,M)
|_{q\to q^i,t\to t^i}}=
\tr_{q,t}(F_U^{\mathrm{bigr}}(\varphi),F_U^{\mathrm{bigr}}(M)). \]
\end{proposition}
\begin{proof}
The proof is virtually identical to that of Proposition \ref{traceprop}; 
\eqref{fundtraceeqn} now becomes the fact that
when $w$ consists of a single
cycle of length $n$ and type $\zeta$,
\begin{equation} \label{gfundtraceeqn} 
\tr_{q,t}(\varphi^n\zeta,M)|_{q\to q^n,t\to t^n}=
\tr_{q,t}(w\varphi^{\otimes n},M^{\otimes n}).
\end{equation}
By the same argument as before, it suffices to prove this when $M$
is one-dimensional, say of bidegree $(a,b)$. Let $x=\tr(\varphi,M)$
and $y=\tr(\zeta,M)$. By our conventions, $M^{\otimes n}$ has bidegree
$(na,nb)$, and $\tr(w,M^{\otimes n})=(-1)^{b(n-1)}y$, since 
if $\tilde{w}\in S_n$ is the $n$-cycle corresponding to $w$,
\[ \varepsilon(\tilde{w},(b,\cdots,b))=(-1)^{b(n-1)}. \]
Thus
\eqref{gfundtraceeqn} asserts that
\[ x^n y (q^n)^a (-t^n)^b= x^n (-1)^{b(n-1)} y q^{na} (-t)^{nb}, \]
which is true.
\end{proof}
\noindent
A special case of this Proposition is that
if $M$ is ungraded (i.e.\ concentrated in bidegree $(0,0)$),
\begin{equation} \label{gtraceeqn}
\ch_{q,t}(U)|_{p_i(\zeta)\to\tr(\varphi^i\zeta,M)}=
\tr_{q,t}(F_U^{\mathrm{bigr}}(\varphi),F_U^{\mathrm{bigr}}(M)).
\end{equation}
The bigrading on $F_U^{\mathrm{bigr}}(M)$ in this case comes purely from that
on $U$.

The sum and product of bigraded $\B_r$-modules are bigraded in the obvious way.
\begin{proposition} \label{gsumprodprop}
If $U,V$ are bounded bigraded $\B_r$-modules, then
\begin{equation*}
\begin{split}
F_{U+V}^{\mathrm{bigr}}\cong F_U^{\mathrm{bigr}}
\oplus F_V^{\mathrm{bigr}},&\ 
F_{U\cdot V}^{\mathrm{bigr}}
\cong F_U^{\mathrm{bigr}}\otimes F_V^{\mathrm{bigr}},\\
\ch_{q,t}(U+V)=\ch_{q,t}(U)+\gch(V),&\ \gch(U\cdot V)=\gch(U)\gch(V).
\end{split}
\end{equation*}
\end{proposition}
\begin{proof}
The isomorphisms of polynomial functors are proved in the same way
as Proposition \ref{functsumprodprop}. The equalities of bigraded
characteristics can either be deduced from these using Proposition
\ref{gtraceprop}, or proved directly.
\end{proof}
\begin{corollary} \label{gsumprodcor}
$\ch_q(U+V)=\ch_q(U)+\ch_q(V), \ch_q(U\cdot V)=\ch_q(U)\ch_q(V)$.
\end{corollary}
\begin{proof}
Simply specialize $t\to 1$ in Proposition \ref{gsumprodprop}.
\end{proof}

If $U$ is a bigraded $\B_r$-module and $V$ a bigraded $\B_1$-module, 
we define the bigraded $\B_r$-module $U\circ V$
by the same formula as before:
\[ (U\circ V)(I)=\bigoplus_{\pi\in\Par(I)_{\B_r}}
\left(U(\pi)\otimes\bigotimes_{\CalO\in\mu_r\setminus\pi}V(\CalO)\right), \] 
but incorporating
the above sign convention on the tensor product 
$\bigotimes_{\CalO\in\mu_r\setminus\pi} V(\CalO)$. That is, if $f:I\isomto I$
is an endomorphism in $\B_r$ preserving $\pi\in\Par(I)_{\B_r}$, then
the action of $(U\circ V)(f)$ on $\bigotimes_{\CalO\in\mu_r\setminus\pi} 
V(\CalO)_{a_\CalO,b_\CalO}$
should have a sign
$\varepsilon(\tilde{w}_f,
(b_{\CalO_1},\cdots,b_{\CalO_{|\mu_r\setminus\pi|}}))$,
where $\CalO_1,\cdots,\CalO_{|\mu_r\setminus\pi|}$ is any ordering of the
$\mu_r$-orbits in $\pi$ and $\tilde{w}_f$ is the permutation of them induced
by $f$.

With $F_V^{(r),\mathrm{bigr}}:\C\mu_r$-$\mathbf{mod}^\mathrm{bigr}\to
\C\mu_r$-$\mathbf{mod}^\mathrm{bigr}$ being the functor induced by 
$F_V^{\mathrm{bigr}}$, we have the analogue of
Theorem \ref{functisomthm}:
\begin{theorem} \label{gfunctisomthm}
If $U$ is a bounded bigraded $\B_r$-module, and $V$ is a bounded
bigraded $\B_1$-module with $V(\emptyset)=\{0\}$, then
$F_{U\circ V}^{\mathrm{bigr}}\cong F_U^{\mathrm{bigr}}\circ 
F_V^{(r),\mathrm{bigr}}$.
\end{theorem}
\begin{proof}
The proof is the same as that of Theorem \ref{functisomthm}; the
reader can check that the signs we have introduced are consistent.
\end{proof}

Finally we come to the bigraded analogue of Theorem 
\ref{plethysmthm}. We define plethysm $\circ:\Aa(r)[q,t]\times\Aa(1)[q,t]_+
\to\Aa(r)[q,t]$ in the same way as in \S1, with the extra rule that
$p_i(\zeta)\circ t = t^i$.
\begin{theorem} \label{gplethysmthm}
If $U$ is a bigraded $\B_r$-module, and $V$ is a bigraded $\B_1$-module with
$V(\emptyset)=\{0\}$, then $\gch(U\circ V)=\gch(U)\circ\gch(V)$.
\end{theorem}
\begin{proof}
As with Theorem \ref{plethysmthm}, it suffices to prove this when
$U$ and $V$ are bounded, and in that case to prove that
the two sides become equal under every specialization of the form
$p_i(\zeta)\to\tr(\varphi^i\zeta,M)$ for $M$ and $\phi$ (ungraded)
as above. Using Proposition \ref{gtraceprop}, \eqref{gtraceeqn}, 
and Theorem \ref{gfunctisomthm}, we have:
\begin{equation*}
\begin{split}
\gch(U\circ V)|_{p_i(\zeta)\to\tr(\varphi^i\zeta,M)}
&=\gtr(F_{U\circ V}^{\mathrm{bigr}}(\varphi),F_{U\circ V}^{\mathrm{bigr}}(M))\\
&=\gtr(F_U^{\mathrm{bigr}}(F_V^{(r),\mathrm{bigr}}(\varphi)),
F_U^{\mathrm{bigr}}(F_V^{(r),\mathrm{bigr}}(M)))\\
&=\gch(U)|_{p_i(\zeta)\to\gtr(F_V^{(r),\mathrm{bigr}}(\varphi)^i\zeta,
F_V^{(r),\mathrm{bigr}}(M))|_{q\to q^i, t\to t^i}}\\
&=\gch(U)|_{p_i(\zeta)\to\gtr(F_V^{\mathrm{bigr}}(\varphi^i\zeta),
F_V^{\mathrm{bigr}}(M))|_{q\to q^i,t\to t^i}}\\
&=\gch(U)|_{p_i(\zeta)\to\gch(V)|_{p_j\to\tr(\varphi^{ij}\zeta^j,M),q\to q^i,
t\to t^i}}\\
&=(\gch(U)\circ\gch(V))|_{p_i(\zeta)\to\tr(\varphi^i\zeta,M)},
\end{split}
\end{equation*}
as required.
\end{proof}
\begin{corollary} \label{gplethysmcor}
With $U$ and $V$ as above, $\ch_q(U\circ V)=\ch_q(U)\circ\ch_q(V)$.
\end{corollary}
\begin{proof}
Clearly the specialization $t\to 1$ respects plethysm.
\end{proof}
%%%%%%%%%%%%%%%%%%%%%%%%%%%%%%%%%%%%%%%%%%%%%%%%%%%%%%%%%%%%%%%%
\section{Proof of Theorems \ref{ginzburgthm} and \ref{mainthm}}
Call a complex variety $X$ \emph{Hodge-even}
if each nonzero cohomology
group $H_c^s(X,\C)$ is a direct sum of pure Hodge structures with even weights.
As mentioned already in \S2, this holds if $X$ is a hyperplane complement
(in which case $H_c^s(X,\C)$ is pure of weight $2s-2\dim X$),
or if $X$ is a nonsingular projective variety with no odd cohomologies
(in which case $H^{2s}(X,\C)$ is pure of weight $2s$). When $X$ is
Hodge-even, we will
regard the total cohomology $H_c^\bullet(X,\C)$ as a bigraded vector space,
where the $t$-grading is given by degree, and the $q$-grading is half the
Hodge weight. 

We say a $\B_r$-variety $X$ is Hodge-even if every $X(I)$ is.
In this case we obtain a bigraded $\B_r$-module
$H_c^\bullet X$ by the rule $H_c^\bullet X(I):=H_c^\bullet(X(I),\C)$.
(Every isomorphism between two varieties induces an isomorphism of their
cohomologies, preserving Hodge weights.)
If all $X(I)$ are projective, we will write $H^\bullet X$ instead
of $H_c^\bullet X$. For instance, if $X$ is actually an $r$-species
(i.e.\ all $X(I)$ are finite), $H^\bullet X$ is just $H^0 X$, placed in
bidegree $(0,0)$. More importantly, we have bigraded $\B_r$-modules
$H_c^\bullet\M(r)$ and $H^\bullet\overline{\M}(r)$.
Clearly we can interpret the series $\CalP(r)$ and $\overline{\CalP}(r)$
defined in \S2 as bigraded characteristics with $t$ set to $1$:
\begin{equation*}
\begin{split}
\CalP(r)&=\ch_q(H_c^\bullet\M(r)),\\
\overline{\CalP}(r)&=\ch_q(H^\bullet\overline{\M}(r)).
\end{split}
\end{equation*}

By the K\"unneth Theorem in mixed Hodge theory,
if $X_1$ and $X_2$ are two Hodge-even varieties, we have an isomorphism 
\[ H_c^\bullet(X_1\times X_2,\C)\cong
H_c^\bullet(X_1,\C)\otimes H_c^\bullet(X_2,\C) \] 
of bigraded vector spaces.
If $X_1=X_2=X$, the tensor product is graded-commutative with respect to the 
degree grading. Thus in
the isomorphism $H_c^\bullet(X^n,\C)\cong H_c^\bullet(X,\C)^{\otimes n}$,
the action of $S_n$ on the left corresponds to the
permutation of factors on the right, together with the sign convention used
in the previous section. As a result, we have the following:
\begin{proposition} \label{hcprop}
\begin{enumerate}
\item Let $X$ and $Y$ be Hodge-even $\B_r$-varieties. Then $X+Y$
and $X\cdot Y$ are Hodge-even, and
\[ H_c^\bullet(X+Y)\cong H_c^\bullet X + H_c^\bullet Y,\
H_c^\bullet(X\cdot Y)\cong H_c^\bullet X \cdot H_c^\bullet Y. \]
\item Let $X$ be a Hodge-even $\B_r$-variety and 
$Y$ a Hodge-even $\B_1$-variety.
Then the $\B_r$-variety $X\circ Y$ is Hodge-even, and
\[ H_c^\bullet(X\circ Y)\cong H_c^\bullet X\circ H_c^\bullet Y. \]
\end{enumerate}
\end{proposition}
\noindent
Here the sum, product, and substitution of bigraded $\B_r$-modules
are defined as in \S6. Moreover:
\begin{proposition} \label{chqprop}
If $X\strat Y$ is a stratification of Hodge-even $\B_r$-varieties,
then $\ch_q(H_c^\bullet X)=\ch_q(H_c^\bullet Y)$ in $\Aa(r)[q]$.
\end{proposition}
\begin{proof}
Using the long exact sequence of cohomology with compact supports, we see
that the disjoint unions
$X(I)=\cup_{\alpha\in A(I)} X_\alpha$ induce equalities
\[ \sum_{i}(-1)^i [H_c^i(X(I),\C)]=\sum_{\alpha\in A(I)}\sum_i
(-1)^i [H_c^i(X_\alpha,\C)] \]
in the Grothendieck group of mixed Hodge structures.
After replacing $[H_c^i(X_\alpha,\C)]$ by $[H_c^i(Y_\alpha,\C)]$, 
the result follows.
\end{proof}
\noindent
Note that it is not in general true that $H_c^\bullet X\cong H_c^\bullet Y$;
one needs to take alternating sum with respect to degree.

Applying Proposition \ref{chqprop} to Theorem \ref{1stratthm},
we get the following equation in $\Aa(1)[q]$:
\begin{equation*}
\overline{\CalP}(1)=
\ch_q(H_c^\bullet\mathbb{T}\M(1)).
\end{equation*}
But by Propositions \ref{1plethprop} and \ref{hcprop}, we have
an isomorphism of bigraded $\B_1$-modules:
\begin{equation*}
H_c^\bullet\mathbb{T}\M(1)\cong
H^0 E(1)_1 + H_c^\bullet \M(1)\circ
H_c^\bullet\mathbb{T}\M(1).
\end{equation*}
Taking $\ch_q$ of both sides, and using Corollaries \ref{gsumprodcor} and 
\ref{gplethysmcor}, we get 
\[ \overline{\CalP}(1)=p_1+\CalP(1)\circ\overline{\CalP}(1), \]
which is Theorem \ref{ginzburgthm}.

If $r\geq 2$, Theorem \ref{rstratthm} and
Proposition \ref{chqprop} imply the following equation in
$\Aa(r)[q]$:
\begin{equation*}
\overline{\CalP}(r)=
\ch_q(H_c^\bullet\mathbb{T}_r(\M(r),\M(1))).
\end{equation*}
But by Propositions \ref{rplethprop} and \ref{hcprop}, we have an isomorphism
of bigraded $\B_r$-modules:
\begin{equation*}
\begin{split}
H_c^\bullet\mathbb{T}_r(\M(r),\M(1))&\cong\\
(H^0 E(r)_0 &+ H_c^\bullet\mathbb{T}_r(\M(r),\M(1)))
\cdot (H_c^\bullet \M(r)\circ H_c^\bullet \mathbb{T}\M(1)).
\end{split}
\end{equation*}
Taking $\ch_q$ of both sides, and using Corollaries \ref{gsumprodcor} and 
\ref{gplethysmcor}, we get 
\[ \overline{\CalP}(r)=(1+\overline{\CalP}(r))
(\CalP(r)\circ\overline{\CalP}(1)), \]
which is Theorem \ref{mainthm}.
%%%%%%%%%%%%%%%%%%%%%%%%%%%%%%%%%%%%%%%%%%%%%%%%%%%%%%%%
\section{Proof of Theorems \ref{1lehrerthm} and \ref{lehrerthm}}
In this section we use the machinery
developed in this paper to give slick proofs of Theorems \ref{1lehrerthm} and
\ref{lehrerthm}. In the case
of Theorem \ref{1lehrerthm} this was done by Getzler in 
\cite[\S5]{mixedhodge}, but as his
paper is not published we will include this case here.
We need a result from the theory of species (see \cite[\S2.5]{bll}),
which says that the plethystic inverse of $Z_{E(1)_+}$ 
in $\Aa(1)$ is
\begin{equation} \label{logeqn} 
\sum_{d\geq 1}\frac{\mu(d)}{d}\log(1+p_d).
\end{equation}
This can be proved easily by direct calculation using the formula
\eqref{expcheqn} for $Z_{E(1)}$. 

Now define a $\B_r$-variety $X(r)$ by the rule
\[ X(r,I):=
\{(z_i)\in\C^I\,|\,z_{\zeta.i}=\zeta z_i,\ \forall i\in I,\zeta\in\mu_r\}. \]
(If $I=\emptyset$, this is interpreted to be $\{0\}$.)
Thus if $r\geq 2$, $X(r)$ coincides with $V(r)$. Clearly $X(r,n)\cong\C^n$.
Define a $\B_r$-variety
$M(r)$ by letting $M(r,I)$ be the complement in
$X(r,I)$ of the hyperplanes with equations $z_i=z_{i'}$ for $i\neq i'$ in $I$.
In particular, $M(r,n)$ can be identified with the complement in $\C^n$
of the hyperplanes $\{z_i=\zeta z_j\}$ for $i\neq j\in [n]$, $\zeta\in\mu_r$,
and (in the $r\geq 2$ case only) $\{z_i=0\}$ for $i\in [n]$.
We define the
bigraded $\B_r$-modules $H_c^\bullet X(r)$, $H_c^\bullet M(r)$
as in the previous section.
\begin{proposition}
In $\Aa(r)[q]$, we have the equation
\[ \ch_q(H_c^\bullet X(r))=
\exp(\sum_{i\geq 1}\frac{q^i}{ri}\sum_{\zeta\in\mu_r}p_i(\zeta)). \]
\end{proposition}
\begin{proof}
Define a $\B_1$-variety $X$ by the rule
\[ X(I)=\left\{\begin{array}{cl}
\C,&\text{ if $|I|=1$,}\\
\emptyset,&\text{ otherwise.}
\end{array}\right. \]
with $X(f)$ being the identity map if $f:I\isomto J$, $|I|=|J|=1$.
Obviously $\ch_q(H_c^\bullet X)=qp_1$. Now by the
definition of substitution,
we have $X(r)\cong E(r)\circ X$, whence
$H_c^\bullet X(r)\cong H^0 E(r)\circ H_c^\bullet X$,
whence $\ch_q(H_c^\bullet X(r))=Z_{E(r)}\circ qp_1$.
Using the formula \eqref{expcheqn} for $Z_{E(r)}$, we get the result.
\end{proof}

It is straightforward to compare the hyperplane complement $M(r,I)$
with the projective hyperplane complement $\M(r,I)$. If $I=\emptyset$, 
$M(r,I)$ consists of a single point whereas $\M(r,I)$ is empty.
If $|I|=1$, $M(1,I)$ is an affine line, whereas $\M(1,I)$
is empty. In all other cases, we have a bundle
projection $M(r,I)\to\M(r,I)$, with fibre $\C\rtimes\C^\times$ if $r=1$
or $\C^\times$ if $r\geq 2$. Hence
\begin{equation} 
\ch_q(H_c^\bullet M(r))=\left\{\begin{array}{cl}
1+qp_1+q(q-1)\CalP(1),&\text{ if $r=1$,}\\
1+(q-1)\CalP(r),&\text{ if $r\geq 2$.}\\
\end{array}\right.
\end{equation}
So Theorem \ref{1lehrerthm} is equivalent to the following.
\begin{theorem} \label{penultthm}
In $\Aa(1)[q]$, we have the equation
\[ \ch_q(H_c^\bullet M(1))=\prod_{n\geq 1}(1+p_n)^{R_n/n}, \]
where $R_n=\sum_{d|n}\mu(d)q^{n/d}$.
\end{theorem}
\begin{proof}
Getzler's proof is as follows. We can stratify $X(1,I)$
according to which of the coordinates of the $I$-tuple are equal.
In our terminology, this gives a stratification $X(1)\strat M(1)\circ E(1)_+$
of $\B_1$-varieties. Applying Propositions \ref{chqprop} and \ref{hcprop} and
Corollary \ref{gplethysmcor}, we get
\begin{equation*}
\begin{split}
\ch_q(H_c^\bullet X(1))&=\ch_q(H_c^\bullet M(1)\circ H^0 E(1)_+)\\
&=\ch_q(H_c^\bullet M(1))\circ Z_{E(1)_+}.
\end{split}
\end{equation*}
Now taking plethysm on the right with the plethystic inverse
of $Z_{E(1)_+}$, given in \eqref{logeqn}, we find that
\begin{equation*}
\begin{split}
\ch_q(H_c^\bullet M(1))&=\ch_q(H_c^\bullet X(1))
\circ(\sum_{d\geq 1}\frac{\mu(d)}{d}\log(1+p_d))\\
&=\exp(\sum_{i\geq 1}\frac{q^i}{i}p_i)
\circ(\sum_{d\geq 1}\frac{\mu(d)}{d}\log(1+p_d))\\
&=\exp(\sum_{i,d\geq 1}\frac{q^i\mu(d)}{id}\log(1+p_{id}))\\
&=\exp(\sum_{n\geq 1}\frac{R_n}{n}\log(1+p_n))\\
&=\prod_{n\geq 1}(1+p_n)^{R_n/n},
\end{split}
\end{equation*}
as claimed.
\end{proof}
\begin{corollary} 
If $w\in S_n$ has $a_i$ cycles of length $i$, then
\[ \sum_{s=n}^{2n}(-1)^s\tr(w,H_c^s(M(1,n),\C))q^{s-n}
=\prod_{i\geq 1}R_i(R_i-i)\cdots(R_i-(a_i-1)i). \]
\end{corollary}
\begin{proof}
By definition of $\ch_q$, the left-hand side is the order of the
centralizer of $w$ multiplied by the coefficient of $p_w$ in
$\ch_q(H_c^\bullet M(1))$. By \eqref{centordereqn}, the order of the
centralizer is $\prod_{i\geq 1}a_i!i^{a_i}$. By
Theorem \ref{penultthm}, the coefficient of $p_w=\prod_{i\geq 1}p_i^{a_i}$ in
$\ch_q(H_c^\bullet M(1))$ is $\prod_{i\geq 1}\binom{R_i/i}{a_i}$.
The result follows.
\end{proof}
\noindent
Compare \cite[Theorem 5.5]{lehrerone}.

To prove Theorem \ref{lehrerthm} we need to
modify this argument only slightly.
\begin{theorem} \label{ultthm}
In $\Aa(r)[q]$ we have the equation
\[ \ch_q(H_c^\bullet M(r))=\prod_{\substack{n\geq 1\\\theta\in\mu_r}}
(1+p_{n}(\theta))^{R_{r,n,\theta}/rn}, \]
where
$R_{r,n,\theta}=\sum_{d|n}
|\{\zeta\in\mu_r\,|\,\zeta^d=\theta\}|\,
\mu(d)(q^{n/d}-1)$.
\end{theorem}
\begin{proof}
We can stratify $X(r,I)$ according to which of the coordinates
are zero and which of the others are equal. This gives a
stratification
\[ X(r)\strat E(r)\cdot (M(r)\circ E(1)_+) \]
of $\B_r$-varieties. Applying Propositions \ref{chqprop} and \ref{hcprop} and
Corollaries \ref{gsumprodcor} and \ref{gplethysmcor}, we get
\begin{equation*}
\begin{split}
\ch_q(H_c^\bullet X(r))&=\ch_q(H^0 E(r)\cdot 
(H_c^\bullet M(r)\circ H^0 E(1)_+))\\
&=Z_{E(r)}
(\ch_q(H_c^\bullet M(r))
\circ Z_{E(1)_+}).
\end{split}
\end{equation*}
So
\begin{equation*}
\begin{split}
\ch_q(H_c^\bullet M(r))&=(\ch_q(H_c^\bullet X(r))
Z_{E(r)}^{-1})
\circ(\sum_{d\geq 1}\frac{\mu(d)}{d}\log(1+p_d))\\
&=\exp(\sum_{i\geq 1}\frac{q^i-1}{ri}\sum_{\zeta\in\mu_r}
p_i(\zeta))\circ(\sum_{d\geq 1}\frac{\mu(d)}{d}\log(1+p_d))\\
&=\exp(\sum_{i,d\geq 1}\sum_{\zeta\in\mu_r}\frac{(q^i-1)\mu(d)}{rid}
\log(1+p_{id}(\zeta^d)))\\
&=\exp(\sum_{n\geq 1}\sum_{\theta\in\mu_r}\frac{R_{r,n,\theta}}{rn}
\log(1+p_n(\theta)))\\
&=\prod_{\substack{n\geq 1\\\theta\in\mu_r}}(1+p_n(\theta))
^{R_{r,n,\theta}/rn},
\end{split}
\end{equation*}
as claimed.
\end{proof}
\begin{corollary}
If $w\in W(r,n)$ is in the conjugacy class corresponding to $(a_i(\zeta))$,
then
\begin{equation*}
\begin{split} 
\sum_{s=n}^{2n}(-1)^s\tr(w,H_c^s(M(r,n),\C))q^{s-n}&=\\
\prod_{\substack{i\geq 1\\\zeta\in\mu_r}}
R_{r,i,\zeta}(R_{r,i,\zeta}-ri)\cdots & (R_{r,i,\zeta}-(a_i(\zeta)-1)ri).
\end{split}
\end{equation*}
\end{corollary}
\begin{proof}
Again, the left-hand side is the order of the centralizer of $w$
multiplied by the coefficient of $p_w$ in $\ch_q(H_c^\bullet M(r))$.
Using \eqref{centordereqn} and Theorem \ref{ultthm}, we get the result.
\end{proof}
\noindent
Compare the results on this group in \cite{lehrertwo}.
%%%%%%%%%%%%%%%%%%%%%%%%%%%%%%%%%%%%%%%%%%%%%%%%%%%%%%%%

%%%%%%%%%%%%%%%%%%%%%%%%%%%%%%%%%%%%%%%%%%%%%%%%%%%%%%%%
\end{document}